%% file: m3-II-1.tex
\input m3-macs

\input xy
\xyoption{all}

\pageno=199

\tinfo II.1.199-213

\SetTFLinebox{\gtp }
\SetFLinebox{\gtv3 }
\SetHLinebox{\issn}

\H 1. Higher dimensional local fields and $L$-functions 

A. N. Parshin

\SetAuthorHead{A. N. Parshin}
\SetTitleHead{Part II. Section 1. Higher  local fields and $L$-functions  \qquad\qquad}

\HH 1.0. Introduction

\HHH 1.0.1

Recall \cite{P1}, \cite{FP} that if $X$ is a scheme of dimension $n$
and
$$X_0\subset X_{1}\subset\dots X_{n-1}\subset X_n=X$$
is a flag of irreducible subschemes ($\dim (X_i)=i$),
then one can define a ring 
$$K_{X_0,\dots,X_{n-1}}$$ associated to the flag.  In the case where
everything is 
regularly embedded, the ring is an $n$-dimensional local field.
Then one can form an adelic object
$$\Bbb A_X= {\prod}' K_{X_0,\dots,X_{n-1}}$$
where the product is taken
over all the flags  with respect to certain restrictions on components
of adeles \cite{P1}, \cite{Be}, \cite{Hu}, \cite{FP}.

\eg Example

Let $X$ be an algebraic projective irreducible surface over a field $k$
and let $P$ be a closed point of $X$, $C\subset X$ be an irreducible
curve such that $P\in C$.

If  $X$ and $C$ are smooth at $P$, then we let $t\in {\Cal O}_{X,P}$ be a local
equation of $C$ at $P$ and
$u\in {\Cal O}_{X,P}$ be such that $u|_C\in {\Cal O}_{C,P}$
is a local parameter at $P$. Denote by $\Cal C$ the ideal  defining the
curve $C$ near $P$. Now we can introduce a two-dimensional local field
$K_{P,C}$ attached
to the pair $P, C$ by the following procedure  including completions and
localizations:
$$
\matrix
{\widehat{\Cal O}}_{X,P} & = & k(P)[[u,t]] \supset \Cal C = (t)\\
\vert &  & \\
({\widehat{\Cal O}}_{X,P})_{\Cal C} & = &
\text{discrete valuation ring with residue field $k(P)((u))$}\\
\vert &  & \\
{\widehat{\Cal O}}_{P,C}:= \widehat{({\widehat{\Cal O}}_{X,P})}_{\Cal C} & = & k(P)((u))[[t]]\\ 
\vert &  & \\
K_{P,C}:=\text{Frac}\,({\widehat{\Cal O}}_{P,C}) & =& k(P)((u))((t)) 
\endmatrix
$$
Note that the left hand side construction is meaningful 
{\it without} any smoothness condition.

\smallskip

Let  $K_P$ be the minimal subring of $K_{P,C}$ 
which contains $k(X)$ and ${\widehat{\Cal O}}_{X, P}$. 
The ring $K_P$ is not a field in general. 
Then $K\subset K_P\subset K_{P,C}$ and  there is another intermediate 
subring $K_C =\text{Frac}\,({\Cal O}_C) \subset K_{P,C}$.
Note that in dimension 2 there is a duality between points $P$ 
and curves $C$ (generalizing the classical duality between points and 
lines in projective geometry). 
 We can compare the structure 
of adelic components in dimension one and two:

$$
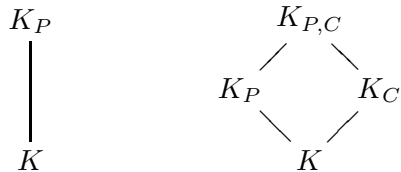
\xymatrix@!0{
K_P\ar@{-}[dd] & & &  & K_{P,C}\ar@{-}[dl] \ar@{-}[dr]  &
\\
 & & & K_P\ar@{-}[dr]  & & K_C \ar@{-}[dl]
\\
K & & & & K &
}
$$

\endeg

\HHH 1.0.2

In the one-dimensional case for every 
character $\chi\colon \Gal(K^{\ab}/K)\to {\Bbb C\,}^*$ we have
the composite 
$$\chi'\colon {\Bbb A}^* = {\prod}'  K_x^* 
@>\text{reciprocity map}>> 
\Gal(K^{\ab}/K) @>{\chi}>> {\Bbb C\,}^* . $$ 
J. Tate \cite{T} and independently K. Iwasawa introduced 
an analytically defined $L$-function 
$$L(s,\chi,f)=\int_{{\Bbb A}^*} f(a)\chi'(a)|a|^s d^*a,$$
where $d^*$ is a Haar measure on ${\Bbb A}^*$ and the function $f$ belongs to
the 
Bruhat--Schwartz space of functions on $\Bbb A$
(for the definition of this space see for instance
\cite{W1, Ch. VII}).  For a special 
choice of $f$ and $\chi = 1$ we get the $\zeta$-function of the scheme $X$
$$
\zeta_X(s) = \prod_{x \in X} (1 - N(x)^{-s})^{-1},
$$
if $\dimm(X) = 1$ (adding the archimedean multipliers if necessary).
Here $x$ runs through the closed points of the scheme $X$
and $N(x) = |k(x)|$. 
The product converges for $\text{Re}(s) > \dimm X$. 
For $L(s,\chi,f)$  they proved the analytical continuation to 
the whole $s$-plane and the functional equation 
$$ 
L(s,\chi,f) = L(1-s,\chi^{-1},{\widehat f}),
$$
using Fourier
transformation ($f \mapsto {\widehat f}$) on the space $\Bbb A_X$ (cf. 
\cite{T}, \cite{W1}, \cite{W2}). 

\HHH 1.0.3

Schemes can be classified according to their dimension
$$
\alignat 3
& \dimm(X)\qquad\qquad & &\text{\rm geometric case} \qquad & &\text{\rm arithmetic case}\\
& \dots \qquad\qquad &&\dots &&\dots \\
& 2\qquad\qquad  &&\text{\rm algebraic surface }/\Bbb F_q \qquad&& \text{\rm arithmetic surface } \\
& 1 \qquad\qquad &&\text{\rm algebraic curve }/\Bbb F_q \qquad &&\text{\rm arithmetic curve}
\\ 
& 0 \qquad\qquad && \Spec(\Bbb F_q)\qquad  && \Spec(\Bbb F_1) 
\endalignat 
$$
where $\Bbb F_1$ is the ``field of one element''.

The analytical method works for the row of the diagram corresponding to dimension one.
The problem to prove analytical continuation 
and functional equation for the $\zeta$-function of arbitrary scheme $X$ 
(Hasse--Weil conjecture)
was formulated by A. Weil \cite{W2} as a generalization of the previous 
Hasse conjecture for algebraic curves over fields of algebraic numbers,  
 see \cite{S1},\cite{S2}.  
It was solved in the geometric situation by A. Grothendieck  
 who employed cohomological methods \cite{G}.  
Up to now there is no extension 
of this method to arithmetic schemes  (see, however, \cite{D}). 
On the other hand, a remarkable property of the Tate--Iwasawa method
is that it can be simultaneously applied  to the fields of algebraic 
numbers (arithmetic situation) and to the algebraic curves over a 
finite field (algebraic situation). 

For a long time the author has been advocating
(see, in particular, \cite{P4}, \cite{FP})
the following:

\th Problem

Extend Tate--Iwasawa's analytic method to higher 
dimensions.
\endth

The higher adeles were introduced exactly for this purpose. In 
dimension one the adelic groups $\Bbb A_X$ and $\Bbb A_X^*$ are locally 
compact groups and thus we can apply the classical harmonic analysis. The 
starting point for 
that is the measure theory on locally compact local fields such as
$K_P$ for the schemes $X$ of dimension 1. So we have the following: 

\th Problem

Develop  a measure theory  and harmonic analysis 
on $n$-dimensional local fields.
\endth

Note that $n$-dimensional local fields  are not locally compact topological spaces for $n>1$ 
and by Weil's theorem the existence of the Haar measure on a topological group 
implies its locally compactness \cite{W3, Appendix 1}. 

In this work  several first steps in answering these problems are described. 

\HH 1.1. Riemann--Hecke method

When one tries to write the $\zeta$-function of a scheme $X$ as a product over
local fields attached to the flags of subvarieties one meets  the following obstacle. 
For dimension greater than one the local fields are parametrized by flags and not by 
the closed points itself as in the Euler product. This problem is primary to 
any problems with the measure and integration. 
I think  we have to return 
to the case of dimension one and reformulate the Tate--Iwasawa method. 
Actually, it means that we have to return to the Riemann--Hecke approach 
\cite{He} known long before the work of Tate and Iwasawa. Of course, it was 
the starting point for their approach. 
 
The main point is a reduction of the integration over ideles to integration
over a single (or finitely many) local field.

Let $C$ be a smooth irreducible complete  curve defined 
over a field $k = \Bbb F_q$.

 Put $K=k(C)$.
For a closed point $x\in C$ denote by $K_x$ the fraction field 
of the completion $\widehat{\Cal O}_x$ of the local ring 
${\Cal O}_x$. 

Let $P$ be a fixed smooth $k$-rational point of $C$. Put $U=C\setminus P$, 
$A=\Gamma(U,\Cal O_C)$. Note that $A$ is a discrete subgroup of $K_P$. 

A classical method to calculate $\zeta$-function is to write it as a 
Dirichlet series instead of the Euler product: 
$$\zeta_C(s)=\sum_{I\in\text{\rm Div\, }({\Cal O}_C)} |I|_C^s$$
where $\text{\rm Div\, }({\Cal O}_C)$ is the semigroup of effective divisors, 
$I=\sum_{x\in X} n_xx$, $n_x\in \Bbb Z$ and $n_x=0$ for almost all $x\in C$,
$$
|I|_C=\prod_{x\in X} q^{-\sum n_x |k(x):k|}. 
$$
Rewrite $\zeta_C(s)$ as 
$$\zeta_U(s)\zeta_P(s)=\biggl( \sum_{I\subset U} |I|_U^s\biggr )
\biggl( \,\sum_{\text{\rm supp}(I)=P} |I|_P^s\biggr ). 
$$ 

Denote $A'=A\setminus \{0\}$. 
For the sake of simplicity assume that $\Pic(U)= (0)$ and introduce $A''$ 
such that $A''\cap k^*=(1)$ and $A' = A''k^*$. 
Then for every $I\subset U$ there is a unique $b\in A''$ such that $I=(b)$.
We also write $|b|_*=|(b)|_*$ for $*=P, U$.
Then 
from the product formula $|b|_C=1$ we get $|b|_U=|b|_P^{-1}$.
Hence
$$
\zeta_C(s)=\biggl(\sum_{b\in A''} |b|_U^s\biggr)
\biggl(\sum_{m\ge 0}  q^{-ms}\biggr)
=\biggl(\sum_{b\in A''} |b|_P^{-s}\biggr)
\int_{a\in K_P^*} |a|_P^s f_{+}(a) d^*a
$$
where in the last equality we have used local Tate's calculation, 
$f_{+}=i^*\delta_{\widehat{\Cal O}_P}$, 
$i\colon K_P^*\to K_P$,
$\delta_{\widehat{\Cal O}_P}$ is the characteristic function of the subgroup 
${\widehat{\Cal O}_P}$, $d^*({\widehat{\Cal O}}^*_P) = 1$.
Therefore
$$
\aligned
\zeta_C(s)&=\sum_{b\in A''}\int_{a\in K_P^*} |ab^{-1}|_P^s f_{+}(a) d^*a
\\
&= \sum_{b\in A''}\int_{c=ab^{-1}}  |c|_P^s f_{+}(bc) d^*c
=\int_{K_P^*} |c|_P^{s} F(c) d^*c,
\endaligned 
$$
where $F(c)=\sum_{b\in A'} f_{+}(bc)$.

Thus, the calculation of $\zeta_C(s)$ is reduced to
integration over the single local field $K_P$. Then we can proceed further
using the Poisson summation formula applied to the function $F$.

This computation can be rewritten in a more functorial way as follows
$$
\zeta_C(s)=\langle|\,\,\,|^s,f_0\rangle_G\cdot \langle|\,\,\,|^s, f_1\rangle_G=
\langle|\,\,\,|^s,i^*(F)\rangle_{G\times G}=\langle|\,\,\,|^s,j_*\circ i^*(F)\rangle_G,
$$
where  $G = K^*_P$, $\langle f,f'\rangle_G=\int_G ff' dg$  and we introduced the functions 
$f_0 = \delta_{A''} =$ sum of Dirac's $\delta_a$ 
over all $a \in A''$ and $f_1 = \delta_{{\Cal O}_P}$ on $K_P$
and the function 
$F = f_0 \otimes f_1$  on $K_P \times K_P$. We also have the 
norm map $|\,\,\,|\colon G \to {\Bbb C\,}^*$, the convolution map $j\colon G \times G \to G$, 
$j(x,y) = x^{-1}y$ and  the inclusion $i\colon G\times G \to K_P \times K_P$. 
 
For the appropriate classes of functions $f_0$ and $f_1$ there are
$\zeta$-functions with a functional equation of the following kind
$$
\zeta(s, f_0, f_1) = \zeta(1-s, \widehat{f_0}, \widehat{f_1}),
$$ 
where ${\widehat f}$ is a Fourier transformation of $f$. We will study the corresponding
spaces of functions and operations like $j_*$ or $i^*$ in subsection 1.3. 

\rk Remark 1

 We assumed that $\Pic(U)$ is trivial. To handle  
the general case one has to consider the curve $C$ with several points removed. 
Finiteness of the $\Pic^0(C)$ implies that we can get an open subset $U$
with this property. 
\endrk

\HH 1.2. Restricted adeles for dimension 2

\HHH 1.2.1

Let us  discuss the situation for dimension one once more. We consider the 
case of the algebraic curve $C$ as above. 

One-dimensional adelic complex 
$$K \oplus \prod_{x\in C}\widehat{\Cal O}_x \to {\prod}'_{x\in C} K_x$$
can be included 
into the following commutative diagram 
$$
\CD
K \oplus \prod_{x\in C}\widehat{\Cal O}_x @>>> {\prod}'_{x\in C} K_x\\
@VVV @VVV \\
K\oplus \widehat{\Cal O}_P @>>> {\prod}'_{x\not=P} K_x/\widehat{\Cal O}_x\oplus K_P
\endCD
$$
where the vertical map induces an isomorphism of cohomologies of the
horizontal complexes. Next, we have a commutative diagram
$$
\CD
K%{\Cal O}_{\eta}
\oplus \widehat{\Cal O}_P @>>> {\prod}'_{x\not=P} K_x/\widehat{\Cal O}_x\oplus K_P\\
@VVV @VVV \\
K/A @>>> {\prod}'_{x\not=P} K_x/\widehat{\Cal O}_x
\endCD
$$ 
where the bottom horizontal arrow is an isomorphism 
(the surjectivity follows from the strong approximation theorem). 
This shows that the complex $A \oplus \widehat{\Cal O}_P \to K_P$ is
quasi-isomorphic to the full adelic complex. 
The construction 
can be extended to an arbitrary locally free sheaf ${\Cal F}$ on $C$
and we obtain that the complex
$$
W \oplus {\widehat{\Cal F}}_P \rightarrow  {\widehat{\Cal F}}_P
\otimes_{{\widehat{\Cal O}}_P} K_P, 
$$ 
where $W = \Gamma({\Cal F}, C\setminus P) \subset K%{\Cal O}_{\eta}
$,
computes the cohomology of the sheaf ${\Cal F}$. 

This fact is essential for the analytical approach to the $\zeta$-function of 
the curve $C$. To understand how to generalize it to higher dimensions
we have to recall another applications of this diagram, in particular, 
the so called Krichever correspondence from the theory of integrable systems. 
 
Let $z$ be a local parameter at $P$, so
$\widehat{\Cal O}_P=k[[z]]$.
The Krichever correspondence assigns points of infinite dimensional
Grassmanians
to $(C,P,z)$ and a torsion free coherent sheaf of ${\Cal O}_C$-modules on
$C$.
In particular, there is an injective map from
classes of triples $(C,P,z)$ to $A\subset k((z))$. In \cite{P5} it was 
generalized to the case of algebraic surfaces using the higher adelic 
language. 

\HHH 1.2.2

Let $X$ be a projective irreducible algebraic surface over a field 
     $k$, $C \subset X$ be an irreducible projective curve, and 
     $P \in C$ be a smooth point on both $C$ and $X$.

In dimension two we  start with  the adelic complex 
$$
\Bbb A_0 \oplus \Bbb A_1 \oplus \Bbb A_2 \rightarrow \Bbb A_{01} \oplus \Bbb A_{02} \oplus \Bbb A_{12}
\rightarrow A_{012}, 
$$ 
where 
$$\aligned
&A_0=K=k(X), \qquad  A_1=\prod_{C\subset X}  \widehat{\Cal O}_C,\qquad  
A_2=\prod_{x\in X}  \widehat{\Cal O}_x,\\
&A_{01}={\prod}'_{C\subset X} K_C,  
A_{02}={\prod}'_{x\in X} K_x, 
A_{12}={\prod}'_{x\in C}  \widehat{\Cal O}_{x,C}, A_{012}=\Bbb A_X={\prod}'     K_{x,C}. 
\endaligned
$$
In fact one can pass to another complex whose cohomologies are the same as of
the adelic complex and which is a generalization of the construction
for dimension one.
 We have to make the following assumptions: 
$P \in C$ is a smooth point on both $C$ and $X$, and the
surface $X \setminus C$ is affine. The desired complex is 
$$ 
A \oplus A_C \oplus \widehat{\Cal O}_P \to B_C \oplus 
B_P \oplus \widehat{\Cal O}_{P,C} \to   K_{P,C} 
$$
where the rings $B_x$, $B_C$, $A_C$ and $A$ have the following meaning. 
Let $x \in C$. Let
\Roster
\Item{} 
     $ B_x = \bigcap\limits_{D \ne C}(K_x \cap
 {\widehat{\Cal O}}_{x,D})$ 
     where the intersection is taken inside $K_{x}$;
\Item{} $ B_C =  K_C \cap  (\bigcap\limits_{x \ne P} B_x)$ 
     where the intersection is taken inside $K_{x,C}$;
\Item{} $ A_C = B_C \cap {\widehat{\Cal O}}_C$, 
  $ A = K \cap (\bigcap_{x \in X\setminus C} {\widehat{\Cal O}}_x)$. 
\endRoster

This can be easily extended to the case of an  arbitrary torsion free coherent
sheaf ${\Cal F}$ on $X$.

\HHH 1.2.3

Returning back to the question about the $\zeta$-function of the surface
$X$ over $k = \Bbb F_q$ we suggest to write it as the product of three Dirichlet
series 
$$
\zeta_X(s) = \zeta_{X\setminus C}(s) \zeta_{C\setminus P}(s) \zeta_P(s) =
\biggl( \sum_{I\subset X\setminus C} |I|_X^s\biggr )
\biggl( \sum_{I\subset C\setminus P} |I|_X^s\biggr )
\biggl( \sum_{I\subset \Spec({\widehat{\Cal O}}_{P,C}) } |I|_X^s\biggr ).
$$ 
Again we can assume that the surface $U = X \setminus C$ has the most 
trivial possible 
structure. Namely, $\Pic(U) = (0)$ and $\text{Ch}(U) = (0)$. Then every rank 2
vector bundle on $U$ is trivial. In the general case one can remove finitely many curves $C$ from 
$X$ to pass to the surface $U$ satisfying these properties (the same idea 
was used in the construction of the higher Bruhat--Tits buildings attached 
to an algebraic surface \cite{P3, sect. 3}). 
 
Therefore any zero-ideal $I$ with support in $X\setminus C$, $C\setminus P$ or $P$ can be defined 
by functions from the rings $A$, $A_C$ and ${\Cal O}_P$, respectively.
The fundamental difference between the case of dimension one and the case of 
surfaces is that zero-cycles $I$ and ideals of finite colength on $X$ are 
not in one-to-one correspondence. 

\rk Remark 2 

In \cite{P2}, \cite{FP} we show that the functional equation 
for the $L$-function on an algebraic surface over a finite field can be 
rewritten using the $K_2$-adeles. Then it has the same shape as the functional 
equation for algebraic curves written in terms of ${\Bbb A}^*$-adeles (as in \cite{W1}). 
\endrk

\HH 1.3. Types for dimension 1

We again discuss the case of dimension one. If $D$ is a divisor on the
curve $C$ then the Riemann--Roch theorem says
$$
l(D) - l(K_C-D) = \degg(D) + \chi({\Cal O}_C),
$$
where as usual $l(D) = \dimm\,\, \Gamma(C, {\Cal O}_X(D))$ and $K_C$ is the canonical
divisor. If $\Bbb A = {\Bbb A\,}_C$ and $\Bbb A_1 = \Bbb A(D)$ then
$$
H^1(C, {\Cal O}_X(D)) = \Bbb A/(\Bbb A(D) + K),
\qquad H^0(C, {\Cal O}_X(D)) = \Bbb A(D) \cap K 
$$
where $K = \Bbb F_q(C)$. 
We have the following topological properties of the groups:
$$
\aligned 
&\Bbb A\\ 
&\Bbb A(D)\\ 
&K\\ 
&\Bbb A(D) \cap K\\\ 
&\Bbb A(D) + K
\endaligned
\qquad
\aligned 
&\text{locally compact group,}\\
&\text{compact group,}\\
&\text{discrete group,}\\
&\text{finite group,}\\
&\text{group of finite index of $\Bbb A$.}
\endaligned 
$$

The group $\Bbb A$ is dual to itself. Fix a rational differential form
$\omega \in \Omega_K^1,~\omega \ne 0$ and an additive character $\psi$ of
$\Bbb F_q$. The following bilinear form
$$
\langle(f_x), (g_x)\rangle = \sum_{x} \res_x(f_xg_x\omega),~(f_x),(g_x) \in \Bbb A
$$ is non-degenerate and defines an auto-duality of $\Bbb A$. 

If we fix a Haar measure $dx$ on $\Bbb A$ then we also have the Fourier transform
$$
f(x) \mapsto  {\widehat f}(x) = \int_{\Bbb A} \psi(\langle x,y\rangle) f(y)dy
$$ 
for functions on $\Bbb A$ and for distributions $F$ defined by the Parseval
equality
$$ ({\widehat F},{\widehat \phi}) = (F, \phi).
$$
One can attach some functions and/or distributions to the subgroups
introduced above
$$ 
\alignat 2 
&\delta_D \qquad &&= \text{the characteristic function of $\Bbb A(D)$} \\
&\delta_{H^1}\qquad &&= \text{the characteristic function of $\Bbb A(D) + K$}\\
&\delta_K \qquad &&= \, \sum_{\gamma \in K} \delta_{\gamma} \quad \text
{where $\delta_{\gamma}$ is the delta-function at the point $\gamma$} \\ 
&\delta_{H^0}\qquad &&= \sum_{\gamma \in \Bbb A(D) \cap K} \delta_{\gamma}.
\endalignat
$$
There are two fundamental rules for the Fourier transform of these
functions
$$
{\widehat \delta}_D = \text{vol}(\Bbb A(D)) \delta_{\Bbb A(D)^{\bot}},
$$
where
$$
\Bbb A(D)^{\bot} = \Bbb A((\omega) - D),
$$
and
$$
{\widehat \delta}_{\Gamma} = \text{vol}(\Bbb A/\Gamma)^{-1}\delta_{\Gamma^{\bot}}
$$
for any discrete co-compact group $\Gamma$.  
In particular, we can apply 
that to $\Gamma = K = \Gamma^{\bot}$. We have 
$$
\aligned
&(\delta_K, \delta_D) = \#(K \cap \Bbb A(D)) = q^{l(D)},\\
&({\widehat\delta}_K, {\widehat\delta}_D) = \text{vol}(\Bbb A(D))\text{vol}(\Bbb A/K)^{-1}
(\delta_K, \delta_{K_C-D}) =q^{\deg D} q^{\chi ({\Cal O}_C)} q^{l(K_C-D)}
\endaligned
$$
and the Parseval equality gives us the Riemann--Roch theorem.

The functions in these computations can be classified according to 
their types. There are four types of functions
which were introduced by F. Bruhat in 1961 \cite{Br}. 

Let $V$ be a finite dimensional vector space over the adelic ring $\Bbb A$ (or
over an one-dimensional local field $K$ with finite residue field $\Bbb F_q$).
We put
$$
\aligned
{\Cal D} &= \{\text{locally constant functions with compact support}\},\\
{\Cal E} &= \{\text{locally constant functions}\},\\
{\Cal D}' &= \{\text{dual to}~{\Cal D}~=~\text{all distributions}\},\\
{\Cal E}' &= \{\text{dual to}~{\Cal E}~=~\text{distributions with compact
support}\}.
\endaligned 
$$
Every $V$ has a filtration $P \supset Q \supset R$ by compact open subgroups
such that all quotients $P/Q$ are finite dimensional vector spaces over
$\Bbb F_q$.
\par\medskip
If $V, V'$ are the vector spaces  over $\Bbb F_q$ of finite dimension then
for every homomorphism
$
i\colon V \rightarrow V'$
there are two maps 
$$
{\Cal F}(V) @>{i_*}>> {\Cal F}(V'),\qquad 
{\Cal F}(V') @>{i^*}>> {\Cal F}(V),
$$
of the spaces ${\Cal F}(V)$ of all functions on $V$ (or $V'$) with values
in $\Bbb C$. Here $i^*$ is the standard inverse image and $i_*$ is defined by
$$
i_*f(v') = \cases
0, &\text{if $v'\notin \im(i)$} \\
\sum_{v\mapsto v'} f(v), &\text{otherwise.}
\endcases 
$$ 
The maps $i_*$ and $i^*$ are dual to each  other. 

We apply these constructions to give a more functorial definition of the
Bruhat spaces. For any triple $P$, $Q$, $R$ as above we have an epimorphism
$i\colon P/R \rightarrow P/Q$
with the corresponding map for functions 
${\Cal F}(P/Q) @>{i^*}>> {\Cal F}(P/R)$ 
and a monomorphism
$j\colon Q/R \rightarrow P/R$ 
with the map for functions 
${\Cal F}(Q/R) @>{j_*}>> {\Cal F}(P/R)$.

Now the Bruhat spaces can be defined as follows 
$$
\aligned 
&{\Cal D} = \inlim_{j_*} \inlim_{i^*} \, {\Cal F}(P/Q), \\
&{\Cal E} = \prlim_{j^*} \inlim_{i^*}  \, {\Cal F}(P/Q), \\
&   {\Cal D'} = \prlim_{j^*} \prlim_{i_*} \, {\Cal F}(P/Q),\\ 
          &{\Cal E'} = \inlim_{j_*}  \prlim_{i_*} \, {\Cal F}(P/Q). 
      \endaligned 
$$ 
The spaces don't depend on the choice of the chain of subspaces 
$P, Q, R$. Clearly we have 
$$
\aligned
&\delta_{D}  \in  {\Cal D(\Bbb A)},\\
&\delta_{K}  \in  {\Cal D'(\Bbb A)}, \\ 
&\delta_{H^0}  \in  {\Cal E'(\Bbb A)}, \\ 
&\delta_{H^1}  \in  {\Cal E(\Bbb A)}.
\endaligned 
$$ 
 We have the same relations for the functions $\delta_{{\Cal O}_P}$ and
$\delta_{A''}$ 
on the group $K_P$ considered in section 1. 

The Fourier transform preserves the spaces ${\Cal D}$ and ${\Cal D'}$ but
interchanges the spaces ${\Cal E}$ and ${\Cal E'}$. Recalling the origin 
of the subgroups from the adelic complex we can say that, in dimension one 
the types of the functions  have the following structure 
$$
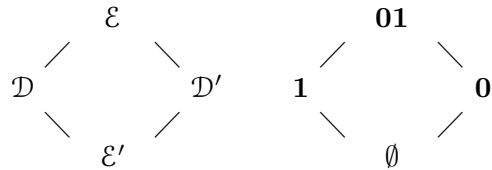
\matrix
&  & \Cal E & & \\
&\diagup &  & \diagdown  \\
\Cal D \!\! \!\!& & & & \!\!\!\! {\Cal D}' \\
&\diagdown  & & \diagup  \\
& & {\Cal E}' & &
\endmatrix
\qquad
\matrix
&  & \boldkey 0\boldkey 1 & & \\
&\diagup &  & \diagdown  \\
\boldkey 1 \!\! \!\!& & & & \!\!\!\! \boldkey 0 \\
&\diagdown  & & \diagup  \\
& & \emptyset & &
\endmatrix
$$
corresponding to the full simplicial division of an edge. The Fourier
transform is a reflection of the diagram  with respect to the middle
horizontal axis.

The main properties of the Fourier transform we need in the proof of the
Riemann-Roch theorem (and of the functional equation of the $\zeta$-function)
can be summarized as the commutativity of the following cube diagram 
\smallskip 
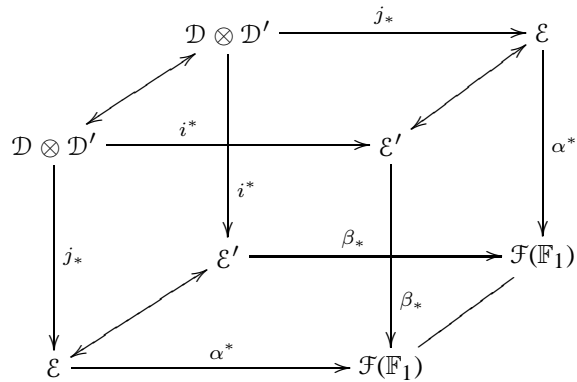
\centerline{\ninepoint%
\xymatrix%@!0
{
&\Cal D\otimes {\Cal D}' \ar@{->}[rr]^{j_*}\ar@{->}[dd]^(0.7){i^*}%|{\vphantom{\int^{\int}}} 
& &\Cal E
\ar@{->}[dd]^{\alpha^*} 
\\
\Cal D\otimes {\Cal D}' \ar@{<->}[ur] 
\ar@{->}[rr]^(0.4){i^*}%|{\hphantom{AAAA}} 
\ar@{->}[dd]^{j_*} 
& & {\Cal E}' \ar@{<->}[ur]
\ar@{->}[dd]^(0.7){\beta_*}%|{\vphantom{\int^{\prod}}} 
\\
& {\Cal E}' \ar@{->}[rr]^(0.4){\beta_*} %|{\hphantom{AAAAAA}}
& & \Cal F(\Bbb F_1)
\\
\Cal E \ar@{->}[rr]^{\alpha^*} 
\ar@{<->}[ur]
&& \Cal F(\Bbb F_1) \ar@{-}[ur]
}}

\noindent coming from the exact sequence 
$$
\Bbb A @>{i}>> \Bbb A \oplus \Bbb A @>{j}>> \Bbb A,
$$
with $i(a) = (a, a),~j(a,b) = a-b$, and the maps
$$
\Bbb F_1 @>{\alpha}>> \Bbb A  @>{\beta}>> \Bbb F_1
$$ 
with $\alpha(0) = 0$, $\beta(a) = 0$.
Here  ${\Bbb F}_{1}$ is the field of one element,
${\Cal F}(\Bbb F_1) =\Bbb C$ and the arrows with  
heads on both ends are the Fourier transforms. 

In particular, the commutativity of the diagram implies
the Parseval equality used above:
$$
\aligned
&\langle \widehat{F}, \widehat{G}\rangle = \beta_*\circ i^*(\widehat{F}\otimes  \widehat{G})\\
&=\beta_*\circ i^*(\widehat{F\otimes G})=
\beta_*\widehat{j_*({F\otimes G})}\\
&=\alpha^*\circ j_*(F\otimes G)=\beta_*\circ i^*(F\otimes G)\\
&=\langle F, G\rangle.
\endaligned
$$

\rk Remark 3

 These constructions can be extended to the function spaces on
the groups $G(\Bbb A)$ or $G(K)$ for a local field $K$ and a group scheme $G$.
\endrk

\HH 1.4. Types for dimension 2

In order to understand the types of functions in the case of dimension 2 we
have to look at the adelic complex of an algebraic surface. We will use
physical notations and denote a space by the discrete index which 
corresponds to it. Thus the adelic complex can be written as 
$$
\emptyset \to {\bold 0} \oplus {\bold 1} \oplus {\bold 2} \to {\bold 0\bold 1} \oplus
{\bold 0\bold 2} \oplus {\bold 1\bold 2} \to {\bold 0\bold 1\bold 2},
$$
where $\emptyset$ stands for the augmentation map corresponding to the inclusion
 of $H^0$. Just as in  the case of dimension one we have a duality of 
$\Bbb A = \Bbb A_{{012}} = {\bold 0\bold 1\bold 2}$ with itself defined by a bilinear form
$$
\langle(f_{x,C}), (g_{x,C})\rangle = \sum_{x,C} \res_{x,C}(f_{x,C}g_{x,C}\,\omega),\qquad
(f_{x,C}),(g_{x,C}) \in \Bbb A 
$$
which is also non-degenerate and defines the autoduality of $\Bbb A$. 

It can be shown that 
$$ 
\Bbb A_{0} = \Bbb A_{01} \cap \Bbb A_{02},\quad 
\Bbb A^{\bot}_{01} = \Bbb A_{01},\quad 
\Bbb A^{\bot}_{02} = \Bbb A_{02},\quad 
\Bbb A^{\bot}_{0} = \Bbb A_{01} \oplus \Bbb A_{02}, 
$$
and so on. 
The proofs depend on the following
residue relations for a rational differential form
$\omega\in \Omega_{k(X)}^2$
$$
\alignat 2
&\text{\rm for all $x\in X$}\qquad && \sum_{C\ni x} \res_{x,C}(\omega)=0,\\
&\text{\rm for all $C\subset X$}\qquad && \sum_{x\in C} \res_{x,C}(\omega)=0. 
\endalignat
$$

We see
that the subgroups appearing in the adelic complex are not
closed under the duality. It means that the set of types in
dimension two will be greater then the set of types coming from the components
of the adelic complex. Namely, we have:

\th Theorem 1 {{\rm (\cite{P4})}} 

Fix a divisor $D$ on an algebraic surface $X$ and 
let $\Bbb A_{12} = \Bbb A(D)$. Consider the lattice ${\Cal L}$ of the
commensurability classes of subspaces in $\Bbb A_X$ generated by
subspaces $\Bbb A_{01}, \Bbb A_{02}, \Bbb A_{12}$.

The lattice ${\Cal L}$ is isomorphic to a free distributive lattice in
three generators and has the  structure shown in the diagram.

$$\ninepoint
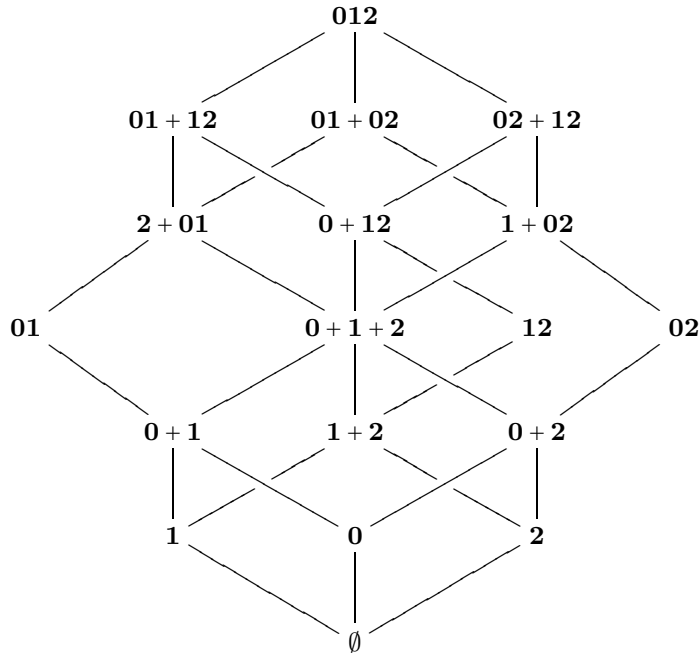
\xymatrix{
& & \bold 0\bold 1\bold 2 \ar@{-}[dl] \ar@{-}[d] \ar@{-}[dr] 
\\
& \bold 0\bold 1+\bold 1\bold 2 \ar@{-}[d] \ar@{-}[dr] & \bold 0\bold 1+\bold 0\bold 2 \ar@{-}[dl]|\hole \ar@{-}[dr]|\hole 
& \bold 0\bold 2+\bold 1\bold 2 \ar@{-}[dl] \ar@{-}[d] 
\\
& \bold 2+\bold 0\bold 1 \ar@{-}[dl] \ar@{-}[dr] & \bold 0+\bold 1\bold 2 \ar@{-}[d] \ar@{-}[dr]|\hole & \bold 1+\bold 0\bold 2 \ar@{-}[dl] \ar@{-}[dr]   
\\
\bold 0\bold 1  \ar@{-}[dr]   &  & \bold 0+\bold 1+\bold 2 \ar@{-}[dl] \ar@{-}[d] \ar@{-}[dr] & \bold 1\bold 2 \ar@{-}[dl]|\hole
& \bold 0\bold 2 \ar@{-}[dl] 
\\
& \bold 0+\bold 1 \ar@{-}[d] \ar@{-}[dr] & \bold 1+\bold 2 \ar@{-}[dl]|\hole \ar@{-}[dr]|\hole &
\bold 0+\bold 2 \ar@{-}[dl] \ar@{-}[d] 
\\
& \bold 1 \ar@{-}[dr]  & \bold 0 \ar@{-}[d] & \bold 2 \ar@{-}[dl] 
\\
& & \emptyset
}
$$

\endth

\rk Remark 4

Two subspaces  $V, V'$ are called commensurable if 
$(V + V')/V \cap V'$ is of  finite dimension.
In the one-dimensional case {\it all} the  subspaces of the 
adelic complex are commensurable (even the subspaces corresponding to
different divisors). In this case we get a free distributive lattice in
two generators (for the theory of lattices see \cite{Bi}).
\endrk

Just as in  the case of curves we can attach to every node some space of
functions (or distributions) on $\Bbb A$. We describe here a particular case
of the construction, namely, the space ${\Cal F}_{02}$ corresponding to
the node ${\bold 0\bold 2}$. Also we will consider not the full adelic group but a single
two-dimensional local field $K = \Bbb F_q((u))((t))$. 

In order to define the space we use the filtration in $K$ by the 
powers $\Cal M^n$ of the maximal ideal $\Cal M = \Bbb F_q((u))[[t]]t$
of $K$ as a discrete valuation (of rank 1) field. 
Then we try 
to use the same procedure as for the local field of dimension 1 (see
above). 

If $P \supset Q \supset R$ are the elements of the filtration
then we need to define the maps 
$$
{\Cal D}(P/R) @>{i_*}>> {\Cal D}(P/Q),\qquad 
{\Cal D}(P/R) @>{j^*}>> {\Cal D}(Q/R)
$$ 
corresponding to an epimorphism
$i\colon P/R \rightarrow P/Q$
and a monomorphism 
$ j\colon Q/R \rightarrow P/R$. 
The map $j^*$ is a restriction of the locally constant functions with
compact support and it is well defined. To define the direct image $i_*$
one needs to integrate along the fibers of the projection $i$. To do that
we have to choose a Haar measure on the fibers for {\it all} $P$, $Q$, $R$
in a consistent way. In other words, we need a system of Haar measures
on all quotients $P/Q$ and by transitivity of the Haar measures in exact
sequences it is enough to do that on all quotients $\Cal M^n/\Cal M^{n+1}$.

Since ${\Cal O}_K/\Cal M = \Bbb F_q((u))=K_1$ we can first choose a Haar measure on the 
residue field $K_1$. It will depend on the choice of a fractional ideal
$\Cal M_{K_1}^i$ normalizing the Haar measure. Next, we have to extend the measure on 
all $\Cal M^n/\Cal M^{n+1}$. Again, it is enough to choose a second local parameter 
$t$ which gives an isomorphism
$$ 
t^n \colon  {\Cal O}_K/\Cal M \to \Cal M^n/\Cal M^{n+1}. 
$$
Having made these choices we can put as above 
$$
{\Cal F}_{02} = \prlim_{j^*} \prlim_{i_*} \,{\Cal D}(P/Q)$$ 
where the space ${\Cal D}$ was introduced in the previous section.

We see that contrary to the one-dimensional case the space ${\Cal F}_{02}$ is not
intrinsically defined. But the choice of all additional data can be easily
controlled.

\th Theorem 2 {{\rm (\cite{P4})}}

The set of the spaces ${\Cal F}_{02}$ is canonically a principal homogeneous
space over the valuation group $\Gamma_K$ of the field $K$.
\endth

Recall that $\Gamma_K$ is non-canonically isomorphic to the lexicographically ordered group $ \Bbb Z \oplus \Bbb Z$. 

One can extend this procedure to other nodes of the diagram of types.
In particular, for ${\bold 0\bold 1\bold 2}$ we get the space which does not depend on 
the choice of the Haar measures. 

The standard subgroup of the type ${\bold 0\bold 2}$ is $B_P = \Bbb F_p[[u]]((t))$ and it is clear 
that 
$$ 
\delta_{B_P} \in {\Cal F}_{02}.
$$
The functions $\delta_{B_C}$ and $\delta_{\widehat{\Cal O}_{P,C}}$ have the
types $\bold 0\bold 1$, $\bold 1\bold 2$ respectively.

\rk Remark 5 

Note that the whole structure of all subspaces in $\Bbb A$
or $K$ corresponding to different divisors or coherent sheaves is more
complicated. The spaces $\Bbb A(D)$ of type ${\bold 1\bold 2}$ are no more  commensurable. 
To describe the whole lattice one has to introduce several equivalence
relations (commensurability up to compact subspace, a locally compact subspace 
and so on).
\endrk

\Bib References

\rf{Be} A. A. Beilinson, { Residues and adeles}, Funct. An. and 
its Appl. {14}(1980), 34--35.

\rf{Bi} G. Birkhoff, {Lattice Theory}, AMS, Providence, 1967.

\rf{Bo}
N. Bourbaki, {Int\'egration}  (Chapitre VII ), Hermann, Paris, 1963.

\rf{Br} 
F. Bruhat,    {Distributions sur un groupe localement compact et 
applications \`a l'\'etude des repr\'esentations des groupes 
$\goth p$-adiques}, Bull. Soc. Math. France {89}(1961), 43--75.

\rf{D} C. Deninger, {Some analogies between number theory and 
dynamical systems}, Doc. Math. J. DMV, Extra volume ICM I (1998), 163--186.

\rf{FP} 
T. Fimmel and A. N. Parshin, { Introduction  into  the  Higher  Adelic 
Theory}, book in preparation.

\rf{G} A. Grothendieck, { Cohomologie $l$-adique and Fonctions
$L$}, S\'eminaire de G\'eom\'etrie Alg\'ebrique du Bois-Marie 1965-66, SGA 5, 
Lecture Notes Math., {589}, Springer-Verlag,
Berlin etc., 1977.

\rf{He}  E. Hecke, { Vorlesungen \"uber die Theorie der
algebraischen Zahlen}, Akad. Verlag, Leipzig, 1923.

\rf{Hu} A. Huber, { On the Parshin-Beilinson adeles for schemes}, 
Abh. Math. Sem. Univ. Hamburg {61}(1991), 249--273.

\rf{P1} A. N. Parshin, { On the arithmetic of two-dimensional schemes.
I. Repartitions and residues}, 
 Izv. Akad. Nauk SSSR Ser. Mat. {40}(1976), 736--773; 
 English translation in Math. USSR Izv. 10 (1976).

\rf{P2} A. N. Parshin, { Chern classes, adeles and L-functions},
J.  reine angew. Math. {341}(1983), 174--192. 

\rf{P3}
A. N. Parshin,  { Higher Bruhat-Tits buildings and vector bundles on an
algebraic surface}, Algebra and Number Theory 
(Proc. Conf.  at Inst. Exp. Math. , Univ. Essen, 1992), 
de Gruyter, Berlin, 1994, 165--192.

\rf{P4} A. N. Parshin, { Lectures at the University Paris-Nord}, 
February--May 1994 and May--June 1996.

\rf{P5} A. N. Parshin, { Krichever correspondence for algebraic 
surfaces}, alg-geom/977911098. 

\rf{S1} J.-P. Serre, { Zeta and $L$-functions}, Proc. Conf. 
Purdue Univ. 1963, New York, 1965,  82--92

\rf{S2} J.-P. Serre, { Facteurs locaux des fonctions zeta 
des vari\'et\'es alg\'ebriques (D\'efinitions et conjectures)}, S\'eminaire 
Delange-Pisot-Poitou, 1969/70, n 19. 

\rf{T} J. Tate, { Fourier analysis in number fields and Hecke's
zeta-functions}, Princeton, 1950, published in Cassels J. W. S., Fr\"olich A., (eds), 
{ Algebraic Number Theory}, Academic Press, London etc., 1967,
305--347. 

\rf{W1} 
A. Weil, {  Basic Number Theory}, Springer-Verlag, 
Berlin etc., 1967.

\rf{W2} A. Weil, { Number theory and algebraic geometry}, Proc.
Intern. Math. Congres, Cambridge, Mass., vol. II,  90--100.

\rf{W3} A. Weil, { L'int\'egration dans les groupes topologiques
et ses applications}, Hermann, Paris, 1940.

\endBib

\Coordinates

Department of  Algebra,
Steklov Mathematical Institute,

Ul. Gubkina 8, Moscow GSP-1, 117966 Russia

E-mail: an\@parshin.mian.su, parshin\@mi.ras.ru
\endCoordinates

\vfill
\pagebreak

\end

%% file: m3-macs.tex
%%%% prevent double loading:
\expandafter\ifx\csname mthreemacsloaded\endcsname\relax\else \fi

\magnification1100
\input amstex

%%% Hack of Plain TeX correction and style macros 
%%% written by Walter Neumann and Larry Siebenmann:

 \catcode`\@=11
 \let\wlog@ld\wlog
 \def\wlog#1{\relax}

 \newif\ifIN@
 \def\m@rker{\m@@rker}
 \def\IN@{\expandafter\INN@\expandafter}
 \long\def\INN@0#1@#2@{\long\def\NI@##1#1##2##3\ENDNI@
    {\ifx\m@rker##2\IN@false\else\IN@true\fi}%
     \expandafter\NI@#2@@#1\m@rker\ENDNI@}
  \newtoks\Initialtoks@  \newtoks\Terminaltoks@
  \def\SPLIT@{\expandafter\SPLITT@\expandafter}
  \def\SPLITT@0#1@#2@{\def\TTILPS@##1#1##2@{%
     \Initialtoks@{##1}\Terminaltoks@{##2}}\expandafter\TTILPS@#2@}
  \newtoks\Trimtoks@

 \def\ForeTrim@{\expandafter\ForeTrim@@\expandafter}
 \def\ForePrim@0 #1@{\Trimtoks@{#1}}
 \def\ForeTrim@@0#1@{\IN@0\m@rker. @\m@rker.#1@%
     \ifIN@\ForePrim@0#1@%
     \else\Trimtoks@\expandafter{#1}\fi}
 
  \def\Trim@0#1@{%
      \ForeTrim@0#1@%
      \IN@0 @\the\Trimtoks@ @%
        \ifIN@
             \SPLIT@0 @\the\Trimtoks@ @\Trimtoks@\Initialtoks@
             \IN@0\the\Terminaltoks@ @ @%
                 \ifIN@
                 \else \Trimtoks@ {FigNameWithSpace}%
                 \fi
        \fi
      }

  %%% Math Bolds
  \font\titlebold=cmbx12 scaled 1200
  \font\twelvebold=cmbx12
  \font\tenbold=cmbx10
  \font\ninebold=cmbx9
  \font\sevenbold=cmbx7
  \font\fivebold=cmbx5

  \input amssym.def \input amssym
  %%% point sizes not loaded by amssym.def:
     \font\titlemsa=msam10 at 14.4pt
     \font\titlemsb=msbm10 at 14.4pt
     \font\titleeufm=eufm10 at 14.4pt
     \font\twelvemsa=msam10 scaled 1200
     \font\twelvemsb=msbm10 scaled 1200
     \font\twelveeufm=eufm10 scaled 1200
     \font\ninemsa=msam9
     \font\ninemsb=msbm9
     \font\nineeufm=eufm9

   %%% Cyrillic fonts (for accents and input, see ams cyr doc)
   \ifx\cyrfam\undefined
   \else
     \immediate\write16{}%
     \message{ !!! cyr fonts already defined. !!! }
     \message{ --- edit out superfluous font defs? }
   \fi
   \newfam\cyrfam
       \font\titlecyr=wncyr10 scaled 1440 %%% no caps?
       \font\twelvecyr=wncyr10 scaled 1200
       \font\tencyr=wncyr10
       \font\ninecyr=wncyr9
       \font\sevencyr=wncyr7
       \font\sixcyr=wncyr6

   %%% Euler script fonts (replacing caligraphic):
   \newfam\eusmfam
       \font\titleeusm=eusm10 scaled 1440
       \font\twelveeusm=eusm10 scaled 1200
       \font\teneusm=eusm10
       \font\nineeusm=eusm9
       \font\seveneusm=eusm7
       
       \font\fiveeusm=eusm5

\let\Cal\cal

 %%% Some fonts not loaded by plain
    \font\ninemrm=cmr9 %% new name for 9 pt math roman
    \font\ninei=cmmi9
    \font\ninesy=cmsy9 
    \skewchar\ninei='177
    \skewchar\ninesy='60

  \font\twelvemrm=cmr10 at 12pt %% new name
  \font\twelvei=cmmi10 at 12pt
  \font\twelvesy=cmsy10 at 12pt
 % \font\twelveex=cmex10 at 12pt

  \font\titlemrm=cmr10 at 14.4pt %% new name
  \font\titlei=cmmi10 at 14.4pt
  \font\titlesy=cmsy10 at 14.4pt
 % \font\titleex=cmex10 at 14.4pt

 %%%% Miscellanious font definitions

  \def\Smallfonts{\ninepoint}

  \def\Hfont{\titlepoint\bf}
  \def\Authorfont{\twelvepoint\it}
  \def\HHfont{\twelvepoint\bf}
  \def\HHHfont{\bf}
  % automatically smaller in 9 point parts
  \def\Bibfont{\tenbf}
  \def\Coordfont{\nineit }% defined in osuPSfnt.sty

  \def \thfont {\bf }
  \def \pffont {\it\itSpacing }
  \def \rkfont {\bf }
  \def \dffont {\bf }
  \def \egfont {\bf }

 %%%%% NINEPOINT %%%%%
 \def\ninepoint{%
  \def\rm{\fam0\ninerm}%
    \textfont0=\ninemrm  \scriptfont0=\sevenrm  \scriptscriptfont0=\fiverm
    \textfont1=\ninei    \scriptfont1=\seveni   \scriptscriptfont1=\fivei
  \def\mit{\fam1\ninei}%
  \def\oldstyle{\fam1\ninei}%
    \textfont2=\ninesy   \scriptfont2=\sevensy  \scriptscriptfont2=\fivesy
    \textfont3=\tenex    \scriptfont3=\tenex    \scriptscriptfont3=\tenex
  \def\it{\fam\itfam\nineit}%
    \textfont\itfam=\nineit
  \def\bf{\ifmmode\fam\bffam\else\ninebf\fi}%
    \textfont\bffam=\ninebold 
    \scriptfont\bffam=\sevenbold 
    \scriptscriptfont\bffam=\fivebold%
  \def\msa{\fam\msafam\ninemsa}%
    \textfont\msafam=\ninemsa 
    \scriptfont\msafam=\sevenmsa
    \scriptscriptfont\msafam=\fivemsa%
  \def\msb{\fam\msbfam\ninemsb}%
    \textfont\msbfam=\ninemsb%
    \scriptfont\msbfam=\sevenmsb%
    \scriptscriptfont\msbfam=\fivemsb%
  \def\eufm{\fam\eufmfam\nineeufm}%
    \textfont\eufmfam=\nineeufm
    \scriptfont\eufmfam=\seveneufm
    \scriptscriptfont\eufmfam=\fiveeufm
   \def\eusm{\fam\eusmfam\nineeusm}%
     \textfont\eusmfam=\nineeusm
     \scriptfont\eusmfam=\seveneusm
     \scriptscriptfont\eusmfam=\fiveeusm
   \def\cyr{\fam\cyrfam\ninecyr}%
     \textfont\cyrfam=\ninecyr
     \scriptfont\cyrfam=\sevencyr
     \scriptscriptfont\cyrfam=\sixcyr%%
  \setbox\strutbox=\hbox{\vrule
      height7pt depth3pt width0pt}%
   \baselineskip=10.8pt\rm}

 \let\eightpoint\ninepoint % we do not use eightpoint

 %%%%% FONTS AT TENPOINT %%%%%
 \def\tenpoint{%
  \def\rm{\fam0\tenrm}%
    \textfont0=\tenmrm \scriptfont0=\sevenrm \scriptscriptfont0=\fiverm%
  \def\mit{\fam1\teni}%
  \def\oldstyle{\fam1\teni}%
    \textfont1=\teni   \scriptfont1=\seveni  \scriptscriptfont1=\fivei%
    \textfont2=\tensy  \scriptfont2=\sevensy \scriptscriptfont2=\fivesy%
    \textfont3=\tenex  \scriptfont3=\tenex   \scriptscriptfont3=\tenex%
  \def\it{\fam\itfam\tenit}%
    \textfont\itfam=\tenit%
  \def\bf{\ifmmode\fam\bffam\else\tenbf\fi}%
    \textfont\bffam=\tenbold% was tenbold for osu
    \scriptfont\bffam=\sevenbold%
    \scriptscriptfont\bffam=\fivebold%
  \def\msa{\fam\msafam\tenmsa}%
    \textfont\msafam=\tenmsa%
    \scriptfont\msafam=\sevenmsa%
    \scriptscriptfont\msafam=\fivemsa%
  \def\msb{\fam\msbfam\tenmsb}%
    \textfont\msbfam=\tenmsb%
    \scriptfont\msbfam=\sevenmsb%
    \scriptscriptfont\msbfam=\fivemsb%
  \def\eufm{\fam\eufmfam\teneufm}%
   \textfont\eufmfam=\teneufm
   \scriptfont\eufmfam=\seveneufm
   \scriptscriptfont\eufmfam=\fiveeufm
   \def\eusm{\fam\eusmfam\teneusm}%
    \textfont\eusmfam=\teneusm
    \scriptfont\eusmfam=\seveneusm
    \scriptscriptfont\eusmfam=\fiveeusm
   \def\cyr{\fam\cyrfam\tencyr}%
    \textfont\cyrfam=\tencyr
    \scriptfont\cyrfam=\sevencyr
    \scriptscriptfont\cyrfam=\sixcyr%%
  \setbox\strutbox=\hbox{\vrule %
      height8.5pt depth3.5ptwidth0pt}%
  \baselineskip=\StdBaselineskip\rm}

 %%%%% FONTS AT TWELVEPOINT %%%%%
 \def\twelvepoint{%
  \def\rm{\fam0\twelverm}%
    \textfont0=\twelvemrm \scriptfont0=\tenmrm \scriptscriptfont0=\sevenrm
    \textfont1=\twelvei   \scriptfont1=\teni   \scriptscriptfont1=\seveni
  \def\mit{\fam1\twelvei}%
  \def\oldstyle{\fam1\twelvei}%
    \textfont2=\twelvesy  \scriptfont2=\tensy  \scriptscriptfont2=\sevensy
    \textfont3=\tenex  \scriptfont3=\tenex  \scriptscriptfont3=\tenex
  \def\it{\fam\itfam\twelveit}%
    \textfont\itfam=\twelveit
  \def\bf{\ifmmode\fam\bffam\else\twelvebf\fi}%
    \textfont\bffam=\twelvebold
    \scriptfont\bffam=\tenbold%
    \scriptscriptfont\bffam=\sevenbold%
  \def\msa{\fam\msafam\twelvemsa}%
    \textfont\msafam=\twelvemsa%
    \scriptfont\msafam=\tenmsa%
    \scriptscriptfont\msafam=\sevenmsa%
  \def\msb{\fam\msbfam\twelvemsb}%
    \textfont\msbfam=\twelvemsb%
    \scriptfont\msbfam=\tenmsb%
    \scriptscriptfont\msbfam=\sevenmsb%
  \def\eufm{\fam\eufmfam\twelveeufm}%
   \textfont\eufmfam=\twelveeufm
   \scriptfont\eufmfam=\teneufm
   \scriptscriptfont\eufmfam=\seveneufm
   \def\eusm{\fam\eusmfam\twelveeusm}%
    \textfont\eusmfam=\twelveeusm
    \scriptfont\eusmfam=\teneusm
    \scriptscriptfont\eusmfam=\seveneusm
   \def\cyr{\fam\cyrfam\tencyr}%
    \textfont\cyrfam=\twelvecyr
    \scriptfont\cyrfam=\tencyr
    \scriptscriptfont\cyrfam=\sevencyr%%
  \setbox\strutbox=\hbox{\vrule
      height10.2pt depth4.55pt width0pt}%
  \baselineskip=14pt\rm}

 %%%%% FONTS AT TITLEPOINT %%%%%
 \def\titlepoint{%
    \textfont0=\titlemrm \scriptfont0=\twelvemrm \scriptscriptfont0=\tenmrm
    \textfont1=\titlei   \scriptfont1=\twelvei   \scriptscriptfont1=\teni
  \def\mit{\fam1\titlei}%
  \def\oldstyle{\fam1\titlei}%
    \textfont2=\titlesy  \scriptfont2=\twelvesy  \scriptscriptfont2=\tensy
    \textfont3=\tenex% math ext not avail in varying sizes??
    \scriptfont3=\tenex
    \scriptscriptfont3=\tenex
  \def\it{\fam\itfam\titleit}%
    \textfont\itfam=\titleit
  \def\bf{\ifmmode\fam\bffam\else\titlebf\fi}%
    \textfont\bffam=\titlebold
    \scriptfont\bffam=\twelvebold%
    \scriptscriptfont\bffam=\tenbold%
  \def\msa{\fam\msafam\titlemsa}%
    \textfont\msafam=\titlemsa%
    \scriptfont\msafam=\twelvemsa%
    \scriptscriptfont\msafam=\tenmsa%
  \def\msb{\fam\msbfam\titlemsb}%
    \textfont\msbfam=\titlemsb%
    \scriptfont\msbfam=\twelvemsb%
    \scriptscriptfont\msbfam=\tenmsb%
  \def\eufm{\fam\eufmfam\titleeufm}%
    \textfont\eufmfam=\titleeufm
    \scriptfont\eufmfam=\twelveeufm
    \scriptscriptfont\eufmfam=\teneufm
   \def\eusm{\fam\eusmfam\titleeusm}%
     \textfont\eusmfam=\titleeusm
     \scriptfont\eusmfam=\twelveeusm
     \scriptscriptfont\eusmfam=\teneusm
   \def\cyr{\fam\cyrfam\tencyr}%
    \textfont\cyrfam=\titlecyr
    \scriptfont\cyrfam=\twelvecyr
    \scriptscriptfont\cyrfam=\tencyr%%
  \setbox\strutbox=\hbox{\vrule
      height12.3pt depth5.54pt width0pt}%
  \baselineskip=16pt\rm}

 %%%% RUNNING HEADINGS
\newbox\AuthorBox\newbox\TitleBox
\newbox\TFLinebox
\newbox\FLinebox
\newbox\HLinebox
\def\SetTFLinebox#1{\setbox\TFLinebox=\hbox{#1}}
\def\SetFLinebox#1{\setbox\FLinebox=\hbox{#1}}
\def\SetHLinebox#1{\setbox\HLinebox=\hbox{#1}}

 \def\SetAuthorHead#1{%
     \setbox\AuthorBox=\hbox{\ninepoint \it 
           \ignorespaces\frenchspacing#1\unskip}}
 \def\SetTitleHead#1{%
     \setbox\TitleBox=\hbox{\ninepoint \it
           \ignorespaces\frenchspacing#1\unskip}}

 %% Italic Spacing Correction
  \def\itSpacing{\relax}
  \def\itSpacingOff{\relax}

  %% Main section headings

 \def\Hrule{\hrule width0pt height0pt}

 %% skip used around proclamations, after section headings,
  % and before subsection-headings:
  \newskip\ProcSkip \ProcSkip 8pt plus2pt minus2pt

 \newskip\LastSkip
 \def\SaveLastSkip{\LastSkip\lastskip}
 \def\RestoreLastSkip{\vskip-\LastSkip\vskip\LastSkip}

 %% Do not indent next paragraph after a header:
 \def\NoindentAfter{\everypar={\setbox0=\lastbox\everypar={}}}

 \long\def\H#1\par#2\par{\notenumber=0 \titlepagetrue%
    {
    \baselineskip=20pt
    \parindent=0pt\parskip=0pt\frenchspacing
    \leftskip=0pt plus .2\hsize minus .3\hsize
    \rightskip=0pt plus .2\hsize minus .3\hsize
 \def\\{\unskip\break}%
    \pretolerance=10000 \Hfont #1\unskip\break
     \vskip7pt\Hrule
\hfill \Authorfont #2\hfill\hfill\unskip}
    \vskip48pt plus 4pt minus 4pt% 60pt=48+12pt
    \par\NoindentAfter\rm}

 \long\def\Hi#1\par#2\par{\notenumber=0 \titlepagetrue%
    {  \baselineskip=0pt  \parindent=0pt\parskip=0pt\frenchspacing
    \leftskip=0pt plus .2\hsize minus .3\hsize
    \rightskip=0pt plus .2\hsize minus .3\hsize
}
    \rm}

 %%% Minor section headings

 \newdimen\PageRemainder
  \def\SetPageRemainder{%\maxdimen case at page tops 12-91 LS
     \PageRemainder=\pagegoal
     \ifdim\PageRemainder=\maxdimen\PageRemainder=\vsize
     \else\advance\PageRemainder by -1\pagetotal\fi}

  \def\Rpt@{}\def\Rpt@@{}

  \long\def\HH#1\par{\par%A
  \SaveLastSkip\removelastskip\goodbreak
  \ifdim\LastSkip<30pt %24pt
     \LastSkip 30pt%24pt 
plus 3pt minus 2pt\fi
  \SetPageRemainder\advance\PageRemainder-\LastSkip
  \ifdim\PageRemainder<150pt
       \edef\Rpt@{remain = \the\PageRemainder\noexpand\\
                pagetotal=\the\pagetotal\noexpand\\
                           pagegoal=\the\pagegoal}%
          \fi
   \ifdim\PageRemainder<65pt %%Head plus 4 lines (approx)
       \ifdim\PageRemainder > 0pt
          \edef\Rpt@@{\noexpand\\
                      Had HH PageRemainder$<$\relax 65pt\noexpand\\
                      Hence forced break!}%
     \vskip 0pt plus .2\PageRemainder\eject %% Pull it out a bit
    \fi\fi
    \vskip\LastSkip\Hrule %%%%%%%%\Hrule added
    \pretolerance=10000\rightskip=0pt plus 3em%B
    \hangafter1 \hangindent=2.2em%
    \noindent
    \HHfont \unskip \Ednote{\Rpt@\Rpt@@}%
            \def\Rpt@{}\def\Rpt@@{}%
            \ignorespaces
            #1\par\rightskip=0pt\pretolerance=\StdPretolerance%
    \NoindentAfter
\tenpoint\rm%
     \medskip \vskip\ProcSkip}%interlineskip adds 2pt to this

  \long\def\HHH#1\par{\par%
  \SaveLastSkip\removelastskip\goodbreak
  \ifdim\LastSkip<\ProcSkip%
     \LastSkip\ProcSkip\fi
  \SetPageRemainder\advance\PageRemainder-\LastSkip
  \ifdim\PageRemainder<150pt
       \edef\Rpt@{remain = \the\PageRemainder\noexpand\\
                pagetotal=\the\pagetotal\noexpand\\
                           pagegoal=\the\pagegoal}%
       \fi
   \ifdim\PageRemainder<48pt  %% 4 lines
        \ifdim\PageRemainder > 0pt
             \edef\Rpt@@{\noexpand\\
                      Had HHH PageRemainder$<$\relax48pt\noexpand\\
                      Hence forced break!}%
       \vskip 0pt plus .2\PageRemainder\eject %% Pull it out a bit
      \fi\fi
   \vskip\LastSkip\par\noindent
   \HHHfont \unskip\Ednote{\Rpt@\Rpt@@}%
  \def\Rpt@{}\def\Rpt@@{}%
  \ignorespaces
   #1\unskip.\quad\rm\ignorespaces
   \ignorepars}

  \long\def\ignorepars#1\par{\def\Test{#1}%
     \ifx\Test\Empty\def\This{\ignorepars}%
        \else\def\This{\Test\par}\fi
           \This}
  \def\Empty{}

 \def\Abstract#1\par{\bgroup\Smallfonts\narrower\HHH #1\par}
 \def\endAbstract{\par\egroup}

 %%%%% Proclamations %%%%%

 \def\ProcBreak{\par%
    \ifdim\lastskip<8pt%
    \removelastskip%
    \penalty-200\vskip\ProcSkip\fi}

 \def\th#1\par{\ProcBreak \noindent
   {\thfont\ignorespaces
    #1\unskip.}\it\itSpacing\kern.4em\ignorepars}%\everymath{\/}

 \def\endth{\ProcBreak\rm\itSpacingOff }%\everymath{}

  %% the theorem statement will be in italic by default

 \def\pf#1\par{\ProcBreak %
    \noindent\pffont#1\unskip.\rm\itSpacingOff{\kern .7em}\ignorepars}

  %% \qed is alternative

  %% A Box for the QED
  \def\qedbox{\hbox{\vbox{
    \hrule width0.2cm height0.2pt
    \hbox to 0.2cm{\vrule height 0.2cm width 0.2pt
             \hfil\vrule height0.2cm width 0.2pt}
    \hrule width0.2cm height 0.2pt}\kern1pt}}

  %% Typing in \qed makes the qedbox right justified:
  \def\qed{\ifmmode\qedbox
    \else\unskip\ \hglue0mm\hfill\qedbox\ProcBreak\fi}

  \def \rk #1\par{\ProcBreak
     \noindent{\rkfont\ignorespaces #1\unskip.}%
     \rm\kern.6em\ignorepars}

  \def \endrk {\medskip\ProcBreak }

  \def \df #1\par{\ProcBreak
     \noindent{\dffont\unskip\ignorespaces #1\unskip.}%
     \rm\kern.6em\ignorepars}

  \def \eg #1\par{\ProcBreak
     \noindent\egfont\unskip\ignorespaces #1\unskip.
     \rm\kern.6em\ignorepars}

  \def \endeg {\medskip\ProcBreak }

  \newdimen\Overhang

   \def\MaxTag@#1#2#3#4#5{\setbox0=\hbox{#4\ignorespaces#2\unskip}%
     \dimen0=\wd0\advance\dimen0 by#3
     \ifdim\dimen0<#5\relax\dimen0=#5\fi
     \expandafter\edef\csname #1Hang\endcsname{\the\dimen0}}

 \def\MaxItemTag#1{\MaxTag@{Item}{#1}{.4em}{\ItemStyle}{\parindent}}%
 \def\MaxItemItemTag#1{%
        \MaxTag@{ItemItem}{#1}{.4em}{\ItemItemStyle}{\parindent}}
 \def\MaxNrTag#1{\MaxTag@{Nr}{#1}{.5em}{\NrStyle}{\parindent}}
 \def\MaxReferenceTag#1{%
        \MaxTag@{Reference}{[#1]}{.6em}{\ninerm}{\parindent}}
 \def\MaxFootTag#1{\MaxTag@{Foot}{#1}{.4em}{\ninerm}{\z@}}

  %% \SetOverhang@ will prevent for tag-text collision
  \def\SetOverhang@{\Overhang=.8\dimen0%
     \advance\Overhang by \wd0\relax%nec!
     \ifdim\Overhang>\hangindent\relax%nec!
       \advance\Overhang by .25\dimen0%
       \Ednote{Tag is pushing text.}\osumess{Tag is pushing text.}%
     \else\Overhang=\hangindent
     \fi}

   %%% \Item
   \def\Item#1{\par\noindent
      \hangafter1\hangindent=\ItemHang
      \setbox0=\hbox{\ItemStyle\ignorespaces#1\unskip}%
      \dimen0=.4em\SetOverhang@% dimen0 is extra space
      \rlap{\box0}\kern\Overhang\ignorespaces}

   %%% \ItemItem
   \def\ItemItem#1{\par\noindent
      \hangafter1\hangindent=\ItemItemHang
      \setbox0=\hbox{\ItemItemStyle\ignorespaces#1\unskip}%
      \dimen0=.4em\SetOverhang@
      \advance\hangindent by \ItemHang
      \kern\ItemHang\rlap{\box0}%
      \kern\Overhang\ignorespaces}

  %%%% \Nr Items without hanging indentation
  \def\Nr#1{\par\noindent\hangindent=\NrHang %not really a hang
    \setbox0=\hbox{\NrStyle\ignorespaces#1\unskip}%
    \dimen0=.5em\SetOverhang@% dimen0 is extra space
    \rlap{\box0}\kern\Overhang
    \hangindent=\z@\ignorespaces}

  %%%% Roster (not compulsory)
  %%  endRoster has to remember \lastskip (e.g. from a \qed) through \egroup.
   \newskip\Rosterskip\Rosterskip 1pt plus1pt %% modifiable
   \def\Roster{\par\ifdim\lastskip<\Rosterskip\removelastskip\vskip\Rosterskip\fi
    \bgroup}
   \def\endRoster{\par\global\edef\LastSkip@{\the\lastskip}\removelastskip
       \egroup\penalty-50\LastSkip\LastSkip@\relax
       \ifdim\LastSkip<\Rosterskip\LastSkip\Rosterskip\fi
       \vskip\LastSkip}%%changed Feb/5/92 WN

 %%%%% Emphasis %%%%%

 %%%%% Vertical spacing %%%%%

 %%%%% References %%%%%

 \def\cite#1{%\relaxnext@
    \def\nextiii@##1,##2\end@{{\frenchspacing\rm 
      \lBr\ignorespaces##1\unskip{\rm,~\ignorespaces##2}\rBr}}%
    \IN@0,@#1@%
    \ifIN@\def\next{\nextiii@#1\end@}\else
    \def\next{{\rm\lBr#1\rBr}}\fi\next}

 %%%%% Bibliography %%%%%

   \def \Bib#1\par{%
       \par\removelastskip\SetPageRemainder
       \ifdim\PageRemainder < 97pt
        \ifdim\PageRemainder > 0pt
        \vfill\eject
       \fi\fi
    \ProcBreak \par\begingroup\parskip=0 pt%
    \goodbreak \vskip 15 pt plus 10 pt
    \noindent\null\hfill\Bibfont% \kern??pt%  (center over what?)
      \ignorespaces #1\unskip\hfill\null\par 
    \frenchspacing \Smallfonts\rm
    \parskip=2.5 pt plus 1 pt minus.5pt%
    \nobreak\vskip 12pt plus 2pt minus2pt\nobreak
    \leftskip=0 pt \baselineskip=10.5pt}

 \def\ReferenceTagSlide{0em}
  \def\ReferenceTagGap{.5em}

  \def \rf#1{\par\noindent
     \hangafter1\hangindent=\ReferenceHang      
     \setbox0=\hbox{\ninerm[\ignorespaces#1\unskip]}%        
     \dimen0=\ReferenceTagGap\SetOverhang@
     \rlap{\kern\ReferenceTagSlide\box0}%       
     \kern\Overhang\ignorespaces}

  \def\ref#1\par#2\par#3\par#4\par{%
     \rf{#1}#2\unskip,\ #3\unskip,\
     #4\unskip.}

  \def\endBib{\par\endgroup\vskip 12pt minus 6pt }

 %%%%% Coordinates %%%%%

  \long\def\Coordinates#1\endCoordinates{%\relax}
 {\par\vskip4pt\def\\{\unskip, }\Coordfont\baselineskip10.5pt\noindent#1}}

 \def\pagecontents{%\TRMargIns new, \Pagetot@l
  \gdef\Pagetot@l{\pagetotal}
  \ifvoid\TRMargIns\else
    \rlap{\kern\hsize\kern10pt\vbox to 0pt{%
         \box\TRMargIns\vss}}\fi
  \ifvoid\topins\else\unvbox\topins\fi
   \dimen@=\dp\@cclv \unvbox\@cclv % open up \box255
   \ifvoid\footins\else % footnote info is present
     \vskip\skip\footins
     \footnoterule
     \unvbox\footins\fi
   \ifr@ggedbottom \kern-\dimen@ \vfil \fi}

  %%%%% Some math accents %%%%%

 \newcount\Ht %pg121; Height register, used in Linefigure & accents

 \def \Acc{\expandafter } %%% What is this for?? WN

 \def\swthat{\raise -1.1 ex\hbox{\sam$\widehat{}$}}
 \def\swttilde{\raise -1.2 ex\hbox{\sam$\widetilde{}$}}
 \def \overdot{{\raise .2 ex \hbox to 0pt {\hss\bf\smash{.}\hss}}}
 \def \overcircle{{\raise .1 ex \hbox to 0pt
    {\sam$\eightpoint\scriptstyle\hss\circ\hss$}}}

 \def \Mathaccent#1#2{{\sam % E.g. #1=\widehat
  \setbox4=\hbox{$\vphantom{#2}$}
  \Ht=\ht4 %pg120
  \setbox5=\hbox{${#1}$}
  \setbox6=\hbox{${#2}$}
  \setbox7=\hbox to .5\wd6{}
  \copy7\kern .1\Ht \raise\Ht sp\hbox{\copy5}\kern-.1\Ht
  \copy7\llap{\box6}
  }}

  \def\SwtCheck #1{
        \ifmmode \check{#1}%
                \else \v {#1}%
                \fi}

 %%  \barpartial : bar over partial is common, tailor!
 \def\barpartial {%
   \kern .17 em
    \overline {\kern -.17 em\partial\kern-.03 em}%
    \kern .03 em}

 %%%   BEtter overline
 
  \def\Overline#1{\setbox1=\hbox{\sam ${#1}$}%
      \ifdim \wd1 > 6pt
    \kern .11 em
    \overline {\kern -.11 em#1\kern-.14 em}
    \kern .14 em
  \else
    \kern .03 em
    \overline {\kern -.03 em#1\kern-.04 em}
    \kern .04 em
  \fi}

 \def\SOverline#1{\setbox1=\hbox{\sam ${#1}$}%
      \ifdim \wd1 > 7pt
    \kern .22 em
    \overline {\kern -.22 em#1\kern-.09 em}%
    \kern .09 em
  \else
    \kern .10 em
    \overline {\kern -.10 em#1\kern-.04 em}%
    \kern .04 em
  \fi}

  %%% Better underline

 \def\Underline#1{\setbox1=\hbox{\sam ${#1}$}%
      \ifdim \wd1 > 6pt
    \kern .11 em
    \underline {\kern -.11 em#1\kern-.14 em}
    \kern .14 em
  \else
    \kern .03 em
    \underline {\kern -.03 em#1\kern-.04 em}
    \kern .04 em
  \fi}

 \def\SUnderline#1{\setbox1=\hbox{\sam ${#1}$}%
      \ifdim \wd1 > 7pt
    \kern .04 em
    \underline {\kern -.04 em#1\kern-.2 em}%
    \kern .2 em
  \else
    \kern .0 em
    \underline {\kern -.0 em#1\kern-.15 em}%
    \kern .15 em
  \fi}

  %%%%% Miscellaneous %%%%%

 \def \Blackbox
   {\leavevmode\hskip .3pt \vbox
   {\hrule height 5pt\hbox{\hskip 4.5pt}}\hskip .5pt}

 \def \XX{\Blackbox\kern.5pt\Blackbox} %% editorial use

  \def\.{.\kern1pt}

  %% unbreakable hyphen (by local change of hyphenchar to -1)
    \def\Hyphen{\edef\this{\the\hyphenchar\font}%
          \hyphenchar\font=-1\char\this\hyphenchar\font=\this}

  %% Prose In Math or Display 
 \ifx\undefined\text
  \def\text#1{\hbox{\rm #1}}\fi %% AMSTeX is more sophisticated

  %% Math Object Names (multi-character math object names)
  %%\nolimits can be cancelled
                                     % by a following \limits if wanted

%%%% Larry's mathsurround hacks:

   \everymath{}  %% initially, but later ...

  \def\PassMath@@{\aftergroup\AfterMath@} %% use \aftergroup LS 5-92

 \let\PassMath@\PassMath@@

 \def\AfterMath@{\futurelet\next\AfterMathMole@}

 \def\AfterMathMole@{%\show\next
      \ifcat\next\space% picks off CR and \par cases too; not \dots
          \def\this{}%{(space)}%
      \else
      \ifcat\next\egroup %
        \def\this{\osumess{Handset mathsurround?? ---(see dollar brace)}}%
      \else
      \def\this{\AAfterMath@}% this minority case slow
      \fi\fi
      \this}

 \def\hyphen@{-}
 \def\paren@{)}
 \def\apostr@{'}

 \def\MSC#1{\kern-.8\mathsurround#1\kern.8\mathsurround}

 \def\AAfterMath@#1{\def\Next{#1}%\show\Next%
    \IN@0\Next @,.;:!?\relax @%
    \ifIN@\def\this{\MSC{\Next}}%
    \else
    \ifx\Next\hyphen@\def\this{\futurelet\next\AfterHyphen@}%
    \else
    \ifx\Next\paren@\def\this{#1}%
    \else 
    \ifx\Next\apostr@\def\this{#1}%
    \else \def\this{\osumess{Handset mathsurround??}%
                 #1}\fi\fi\fi\fi
    \this}

 \def\AfterHyphen@#1{\def\Next{#1}%
   \ifx\Next\hyphen@\def\this{--}\else
   \ifcat\next\space%
   \def\this{\kern-\mathsurround\kern.05em- \Next}\else
   \def\this{\kern-\mathsurround\kern.05em\Hyphen\Next}\fi\fi\this}

%%%% switches
 \def\sam{\mathsurround=\z@\let\PassMath@\relax}  %
 \def\mas{\mathsurround=\StdMathsurround\let\PassMath@\PassMath@@}
 
 \def\Mas{\mathsurround=\StdMathsurround
                \everymath{\PassMath@}\let\PassMath@\PassMath@@}

 \def\m@th{\mathsurround=\z@\everymath{}}%% good general measure

 \def\m@@th{\mathsurround=\z@\everymath={}\let\m@th\relax}

\def\underbar#1{$\setbox\z@\hbox{#1}\dp\z@\z@
      \m@th \underline{\box\z@}$\relax}

\def\mathhexbox#1#2#3{\leavevmode
  \hbox{\m@@th$\m@th \mathchar"#1#2#3$}}

\def\dots{\relax\ifmmode\ldots\else$\m@th\ldots\,$\relax\fi}
   %%% this first \relax is ONLY original

\def\dotfill{\cleaders\hbox{\m@@th$\m@th \mkern1.5mu.\mkern1.5mu$}\hfill}
\def\rightarrowfill{$\m@th\mathord-\mkern-6mu%
  \cleaders\hbox{\m@@th$\mkern-2mu\mathord-\mkern-2mu$}\hfill
  \mkern-6mu\mathord\rightarrow$\relax}
\def\leftarrowfill{$\m@th\mathord\leftarrow\mkern-6mu%
  \cleaders\hbox{\m@@th$\mkern-2mu\mathord-\mkern-2mu$}\hfill
  \mkern-6mu\mathord-$\relax}

\def\downbracefill{$\m@th\braceld\leaders\vrule\hfill\braceru
  \bracelu\leaders\vrule\hfill\bracerd$\relax}
\def\upbracefill{$\m@th\bracelu\leaders\vrule\hfill\bracerd
  \braceld\leaders\vrule\hfill\braceru$\relax}

\def\angle{{\vbox{\m@@th\ialign{$\m@th\scriptstyle##$\crcr
      \not\mathrel{\mkern14mu}\crcr
      \noalign{\nointerlineskip}
      \mkern2.5mu\leaders\hrule height.34pt\hfill\mkern2.5mu\crcr}}}}

\def\big#1{{\m@@th\hbox{$\left#1\vbox to8.5\p@{}\right.\n@space$}}}
\def\Big#1{{\m@@th\hbox{$\left#1\vbox to11.5\p@{}\right.\n@space$}}}
\def\bigg#1{{\m@@th\hbox{$\left#1\vbox to14.5\p@{}\right.\n@space$}}}
\def\Bigg#1{{\m@@th\hbox{$\left#1\vbox to17.5\p@{}\right.\n@space$}}}
\def\n@space{\nulldelimiterspace\z@ \m@th}

\def\root#1\of{\setbox\rootbox\hbox{\m@@th$\m@th\scriptscriptstyle{#1}$}
  \mathpalette\r@@t}
\def\r@@t#1#2{\setbox\z@\hbox{\m@@th$\m@th#1\sqrt{#2}$\relax}
  \dimen@\ht\z@ \advance\dimen@-\dp\z@
  \mkern5mu\raise.6\dimen@\copy\rootbox \mkern-10mu \box\z@}

\def\mathph@nt#1#2{\setbox\z@\hbox{\m@@th$\m@th#1{#2}$}\finph@nt}

\def\mathsm@sh#1#2{\setbox\z@\hbox{\m@@th$\m@th#1{#2}$}\finsm@sh}

\def\@vereq#1#2{\lower.5\p@\vbox{\m@@th\baselineskip\z@skip\lineskip-.5\p@
    \ialign{$\m@th#1\hfil##\hfil$\crcr#2\crcr=\crcr}}}

\def\mathpalette#1#2{\sam\mathchoice{#1\displaystyle{#2}}%
  {#1\textstyle{#2}}{#1\scriptstyle{#2}}{#1\scriptscriptstyle{#2}}\mas}

\def\widehat#1{\setbox\z@\hbox{\sam$#1$}%
 \ifdim\wd\z@>\tw@ em\mathaccent"0\msbfam@5B{#1}%
 \else\mathaccent"0362{#1}\fi}
\def\widetilde#1{\setbox\z@\hbox{\sam$#1$}%
 \ifdim\wd\z@>\tw@ em\mathaccent"0\msbfam@5D{#1}%
 \else\mathaccent"0365{#1}\fi}

 \def\dots{\relax{}
  \ifmmode\def\thedots{\mdots@}\else\def\thedots{\tdots@}\fi %
  \thedots}

 %% \eqno and \leqno need protection
 \let\@ldeqno\eqno\let\@ldleqno\leqno
 \def\eqno{\everymath{}\@ldeqno} \def\leqno{\everymath{}\@ldleqno}

  \let\@ldeqalignno\eqalignno
  \def\eqalignno#1{\sam\@ldeqalignno{#1}\mas}
  \let\@ldeqalign\eqalign
  \def\eqalign#1{\sam\@ldeqalign{#1}\mas}

 \def\overrightarrow#1{\vbox{\m@th\ialign{##\crcr
      \rightarrowfill\crcr\noalign{\kern-\p@\nointerlineskip}
      $\hfil\displaystyle{#1}\hfil$\crcr}}}
 \def\overleftarrow#1{\vbox{\m@th\ialign{##\crcr
      \leftarrowfill\crcr\noalign{\kern-\p@\nointerlineskip}
      $\hfil\displaystyle{#1}\hfil$\crcr}}}
 \def\overbrace#1{\mathop{\vbox{\m@th\ialign{##\crcr\noalign{\kern3\p@}
      \downbracefill\crcr\noalign{\kern3\p@\nointerlineskip}
      $\hfil\displaystyle{#1}\hfil$\crcr}}}\limits}
 \def\underbrace#1{\mathop{\vtop{\m@th\ialign{##\crcr
      $\hfil\displaystyle{#1}\hfil$\crcr\noalign{\kern3\p@\nointerlineskip}
      \upbracefill\crcr\noalign{\kern3\p@}}}}\limits}

  \let\@ldmatrix\matrix
  \let\end@ldmatrix\endmatrix
  \def\matrix{\sam\@ldmatrix}
  \def\endmatrix{\end@ldmatrix\mas}
  \let\@ldgather\gather
  \let\end@ldgather\endgather
  \def\gather{\sam\@ldgather}
  \def\endgather{\end@ldgather\mas}
  \let\@ldalign\align
  \let\end@ldalign\endalign
  \def\align{\sam\@ldalign}
  \def\endalign{\end@ldalign\mas}
  \let\@ldaligned\aligned
  \let\end@ldaligned\endaligned
  \def\aligned{\sam\@ldaligned}
  \def\endaligned{\end@ldaligned\mas}
  \let\@ldtag\tag
  \def\tag{\sam\@ldtag}
   %
  %%% Commutative diagrams : use LamsCD too?

   \let\MinCDArrowWidth\minCDaw@

  %% will be redefined by BoxedEPS.tex

  %%%%% \FigureTitle %%%%%

%%%% End of Larry's mathsurround stuff
%%%% Start of Walter's insert corrections

\newskip\insertskipamount\newskip\inserthardskipamount
\insertskipamount 6pt plus2pt %This is medskipamount without shrink
\inserthardskipamount 6pt
\def\insertskip{\vskip\insertskipamount}
\newcount\SplitTest%        will be set to -1 if a topinsert has split
\def\SetSplitTest{\SplitTest\insertpenalties
  \insert\topins{\floatingpenalty1}%
  \advance\SplitTest-\insertpenalties}
\def\midinsert{\par
 \SaveLastSkip\penalty-150\SetSplitTest\RestoreLastSkip
 \ifnum\SplitTest=-1
  \@midfalse\p@gefalse\else\@midtrue\fi\@ins}
\def\@ins{\par\begingroup\setbox\z@\vbox\bgroup%
  \vglue\inserthardskipamount}
\def\endinsert{\egroup % finish the \vbox
  \if@mid \dimen@\ht\z@ \advance\dimen@\dp\z@
    \advance\dimen@\insertskipamount%            was 12pt (wn)
    \advance\dimen@\pagetotal\advance\dimen@-\pageshrink
    \ifdim\dimen@>\pagegoal\@midfalse\p@gefalse\fi\fi
  \if@mid%
    \ifdim\lastskip<\insertskipamount\removelastskip\insertskip\fi
    \nointerlineskip\box\z@\penalty-200\insertskip
  \else%
    \SaveLastSkip%                                  added (wn)
    \insert\topins{\penalty100 % floating insertion
    \splittopskip\z@skip
    \splitmaxdepth\maxdimen \floatingpenalty\z@
    \ifp@ge \dimen@\dp\z@
    \vbox to\vsize{\unvbox\z@\kern-\dimen@}% depth is zero
    \else \box\z@\nobreak\insertskip\fi}% was \bigskip\fi (wn)
    \RestoreLastSkip%                               added (wn)
   \fi\endgroup}
%% End Walter's insert stuff

 %%%%% Footnotes %%%%%

  \newcount\notenumber
  
  \def\note{\advance\notenumber by 1
    \footnote{\the\notenumber)}}

  \newbox\footbox

 %% The following modifies Plain TeX definitions, qv
  \def\footnote#1{\let\@sf\empty
    %{(the text)} is read later
    \ifhmode\edef\@sf{\spacefactor\the\spacefactor}\/\fi
    \sam${}^{\fam0 #1}$\@sf\vfootnote{#1}}%

  \def\vfootnote#1{\insert\footins\bgroup
     \interlinepenalty100 \splittopskip=1pt
     \floatingpenalty=20000
     \leftskip=0pt\rightskip=0pt%
     \parindent=.3em%% adjust
     \Smallfonts\rm%%osudeG added \Smallfonts
     \FootItem@{#1}%\strut% not nec
     \futurelet\next\fo@t}

  \def\FootItem@#1{\par\hangafter1\hangindent=\FootHang
     \setbox0=\hbox{\ignorespaces#1\unskip}%
     \dimen0=.4em\SetOverhang@% dimen0 is extra space
     \noindent\rlap{\box0}\kern\Overhang\ignorespaces}

  %\MaxFootTag{2)}%% in param file

  \def\fo@t{\ifcat\bgroup\noexpand\next \let\next\f@@t
    \else\let\next\f@t\fi \next}
  \def\f@@t{\bgroup\aftergroup\@foot\let\next}
  \def\f@t#1{\baselineskip=10pt\lineskip=1pt
            \lineskiplimit=0pt #1\@foot}%
     %%osudeG added \baselineskip=? pt\lineskiplimit=0pt
  \def\@foot{%%% special strut osu for end of each note
        \hbox{\vrule height0pt depth5pt width0pt}
        \egroup}
  \skip\footins=12 pt plus 0pt minus 0pt %% was \bigskipamount
    %% space added when footnote is present
  \count\footins=1000 % footnote magnification factor (1 to 1)
  \dimen\footins=8in % maximum footnotes per page

 %%%% Altenatives

  %%  Editorial stuff (delete??)

 \def\osumess#1{\EdSpider{\immediate\write16{Line \the\inputlineno: #1}}}%
 \def\HideEdStuff{\gdef\EdSpider##1{}}

 \font\BigSym=cmmi10 scaled \magstep 4

 \def\change{\InLMargin{\hbox{\BigSym \char63\kern10pt}}}

 \def\beginchange{\InLMargin{\hbox{\sam\twelvepoint$\heartsuit$\kern10pt}}}

 \def\endchange{\InLMargin{\hbox{\sam\twelvepoint$\spadesuit$\kern10pt}}}

 \def\InLMargin#1{\strut\vadjust{%
     \kern-\strutdepth
     \vtop to \strutdepth{%
         \baselineskip\strutdepth
         \llap{\sam$\smash{\hbox{\EdSpider{#1}}}$}\null}}}

 \def\strutdepth{\dp\strutbox}
 \def\strutheight{\ht\strutbox}

 \def\NoteInRMargin#1{\strut\vadjust{%
     \kern-1.001\strutdepth
     \vtop to \strutdepth{%
       \baselineskip\strutdepth
       \vss\rlap{\ninepoint\unskip\hskip\hsize
         \vtop to 0pt{%
           \hsize=16em\hfuzz=\hsize
           \leftskip=10pt%
           \rightskip=0pt plus 10000pt%
           \baselineskip=9.8pt\lineskip=.2pt%
           \let\\\break
           \noindent\EdSpider{#1}\vss}%
                \kern10pt}\hbox{}}%%\hbox{}=\null crucial!!
       }}

 \def\ednote#1{\NoteInRMargin{\tentt #1}}

 \def\cbar{\InLMargin{%
      \dimen0=\strutdepth\advance\dimen0 by \lineskip
      \vrule width 3pt
      height \strutheight depth \dimen0 \kern
      3pt}}

 \def\ccbar{\InLMargin{%
      \dimen0=2\strutdepth\advance\dimen0 by 2\lineskip
      \vrule width 3pt
        height 3\strutheight depth \dimen0 \kern
      3pt}}

 \newinsert\TRMargIns
 \dimen\TRMargIns=\maxdimen
 %\count\TRMargIns=0
 %\skip\TRMargIns=0pt

  \def\Ednote#1{\insert\TRMargIns{%
       \vbox to 0pt{\hsize=140pt\hfuzz=\hsize
           \leftskip=6pt%
           \rightskip=0pt plus 10000pt%
           \baselineskip=9.8pt\lineskip=.2pt%
           \let\\\break
           %\vglue\pagetotal% misplaces notes if inserts are present
           \SetPageRemainder% This ...
           \vglue540pt\vglue-\PageRemainder%  .. is a fix (WN)
           \noindent\EdSpider{\tentt #1}\vss}%
       \smallskip}}

 \def\KillEdStuff{\def\ednote##1{}\def\Ednote##1{}%
      \let\change\relax\let\beginchange\relax\let\endchange\relax
       \let\cbar\relax\let\ccbar\relax}

 %%% Compatibility with osumrip.sty
  %%

 %%% Parameters
  \topskip=12pt
  \newskip\StdBaselineskip % to set \baselineskip
  \StdBaselineskip 12pt
  \lineskip=1.1pt
  \lineskiplimit=.8pt
  \widowpenalty=10000 % 8000 to 10000
  \clubpenalty=10000  % 8000 to 10000
  \abovedisplayskip=6pt plus 1pt minus 1pt
  \abovedisplayshortskip=3pt plus 1.5pt
  \belowdisplayskip=6pt plus 1pt minus 1pt
  \belowdisplayshortskip=5pt plus 1pt minus 1pt
  \hfuzz=1.5pt   % Enable overfull box warnings at console

  \def\StdPretolerance{100}
  \tolerance=\StdPretolerance

  \newdimen\StdMathsurround
  \StdMathsurround=1.5pt % 1pt usual without \Mas
  \mathsurround=\StdMathsurround
  \Mas                   %% sophisticated mathsurround on
 % \Sam                   %% sophisticated mathsurround off

%% marker before English punctuation in displayed math
   \def\prose{\relax\hbox{\kern.6\StdMathsurround}}
  
  \def\StdParskip{0pt}    %% Larry wants {2pt plus 1pt}
  \parskip=\StdParskip
  \parindent=0.5cm
 
%%%% load Times for main body font

  \def\Times{ptmr  } 
  \def\TimesI{ptmri  } 
  \def\TimesB{ptmb  }
  \def\TimesBI{ptmbi  }
  \def\HelveticaN{phvrrn }

  =\Times at 10bp% roman text
  =\TimesB at 10bp% boldface extended
   % slanted roman
  \font\tenit=\TimesI at 10bp% text italic
  =\TimesBI at 10bp

  \font\tenmrm=cmr10  %%new name for math role at full size

%%%%% Fonts at ninepoint %%%%%

    =\Times at 9bp 
    \font\nineit=\TimesI at 9bp 
    =\TimesB at 9bp 
    =\TimesBI at 9bp 

    =\HelveticaN at 9bp 
       % see below

%%%%% Fonts at twelvepoint %%%%%

  =\Times at 12bp
  \font\twelveit=\TimesI at 12bp
  =\TimesB at 12bp

%%%%% Fonts at titlepoint %%%%%

  \font\titleit=\TimesI at 14.4bp
  =\TimesB at 14.4bp

 \SetAuthorHead{AuthorHead} % needs \ninepoint since box set
 \SetTitleHead{TitleHead}  % notably \HeaderFont

%%%% Char adjustments %%%%

  \def\lBr{\raise.125ex\hbox{[\kern.1125ex}}
  \def\rBr{\raise.125ex\hbox{\kern.1125ex]}}

 \setbox\footbox=\hbox{\Smallfonts 2)~}

%% Some optional font dimension and spacing 
%% adjustments beyond this point

%% Correct the lousy spacing of italic f (a hack).

  \bgroup
  \catcode`\@=11 %localised
  \gdef\itSpacing{%
     \xspaceskip=.31em plus.1em minus.05em \sfcode `f=2001
     \itWarning@\let\itWarning@\itWarning@@}
  \gdef\itSpacingOff{%
     \xspaceskip=0pt \sfcode `f=1000
     \let\itWarning@\relax}
   \global\let\itWarning@\relax
  \gdef\itWarning@@{\errmessage{%
  Special italic spacing already in force
  (you have probably omitted an ``endth'').
  See itSpacing macro in osuPSfnt.sty
         }}
  \egroup

 %%% Provisional fontdimen settings
  %%
 \fontdimen1\titlebf=0.0pt
 \fontdimen2\titlebf=3.6135pt
 \fontdimen3\titlebf=2.8908pt
 \fontdimen4\titlebf=1.44539pt
 \fontdimen5\titlebf=6.64882pt
 \fontdimen6\titlebf=14.45398pt
 \fontdimen7\titlebf=1.60439pt

 \fontdimen1\tenbi=0.26794pt
 \fontdimen2\tenbi=2.50937pt
 \fontdimen3\tenbi=2.00749pt
 \fontdimen4\tenbi=1.00374pt
 \fontdimen5\tenbi=4.59717pt
 \fontdimen6\tenbi=10.03749pt
 \fontdimen7\tenbi=1.11415pt

 \fontdimen1\twelverm=0.0pt
 \fontdimen2\twelverm=3.01125pt
 \fontdimen3\twelverm=2.409pt
 \fontdimen4\twelverm=1.2045pt
 \fontdimen5\twelverm=5.39615pt
 \fontdimen6\twelverm=12.045pt
 \fontdimen7\twelverm=1.33699pt

 \fontdimen1\twelveit=0.27731pt
 \fontdimen2\twelveit=3.01125pt
 \fontdimen3\twelveit=2.409pt
 \fontdimen4\twelveit=1.2045pt
 \fontdimen5\twelveit=5.37207pt
 \fontdimen6\twelveit=12.045pt
 \fontdimen7\twelveit=1.33699pt

 \fontdimen1\twelvebf=0.0pt
 \fontdimen2\twelvebf=3.01125pt
 \fontdimen3\twelvebf=2.409pt
 \fontdimen4\twelvebf=1.2045pt
 \fontdimen5\twelvebf=5.5407pt
 \fontdimen6\twelvebf=12.045pt
 \fontdimen7\twelvebf=1.33699pt

 \fontdimen1\tenrm=0.0pt
 \fontdimen2\tenrm=2.50937pt
 \fontdimen3\tenrm=2.00749pt
 \fontdimen4\tenrm=1.00374pt
 \fontdimen5\tenrm=4.49678pt
 \fontdimen6\tenrm=10.03749pt
 \fontdimen7\tenrm=1.11415pt

 \fontdimen1\tenit=0.27731pt
 \fontdimen2\tenit=2.50937pt
 \fontdimen3\tenit=2.00749pt
 \fontdimen4\tenit=1.00374pt
 \fontdimen5\tenit=4.47672pt
 \fontdimen6\tenit=10.03749pt
 \fontdimen7\tenit=1.11415pt

 \fontdimen1\tenbf=0.0pt
 \fontdimen2\tenbf=2.50937pt
 \fontdimen3\tenbf=2.00749pt
 \fontdimen4\tenbf=1.00374pt
 \fontdimen5\tenbf=4.61723pt
 \fontdimen6\tenbf=10.03749pt
 \fontdimen7\tenbf=1.11415pt

 \fontdimen1\ninerm=0.0pt
 \fontdimen2\ninerm=2.25842pt
 \fontdimen3\ninerm=1.80673pt
 \fontdimen4\ninerm=0.90337pt
 \fontdimen5\ninerm=4.0471pt
 \fontdimen6\ninerm=9.03374pt
 \fontdimen7\ninerm=1.00273pt

 \fontdimen1\nineit=0.27731pt
 \fontdimen2\nineit=2.25842pt
 \fontdimen3\nineit=1.80673pt
 \fontdimen4\nineit=0.90337pt
 \fontdimen5\nineit=4.02904pt
 \fontdimen6\nineit=9.03374pt
 \fontdimen7\nineit=1.00273pt

 \fontdimen1\ninebf=0.0pt
 \fontdimen2\ninebf=2.25842pt
 \fontdimen3\ninebf=1.80673pt
 \fontdimen4\ninebf=0.90337pt
 \fontdimen5\ninebf=4.15552pt
 \fontdimen6\ninebf=9.03374pt
 \fontdimen7\ninebf=1.00273pt

 %%% \SetExtraSpaces \MaxSpaceFactor \SetSpaceFactors
  %%  See TeXbook, page 76.

 \newcount\MaxSpaceFactor
 \MaxSpaceFactor=3000 %% to reset later

 %%%%% Tag styles and (hang-) indents
 \def\ItemStyle{\rm}
 \def\NrStyle{\rm}
 \def\ItemItemStyle{\rm}

 %% Analog dimensioning, convenient for local modifications:
 \MaxItemTag{(iii)}
 \MaxItemItemTag{(iii)}
 \MaxNrTag{(2)}
 \MaxFootTag{2)}
 % \MaxReferenceTag{AaaAA} % for biblio
 \def\ReferenceHang{30pt}

 \catcode`\@=\active

%%%%% End of hack of Neumann-Siebenmann macros

\loadbold

=\Times  
=\Times scaled750
=\Times scaled650
\font\rms=\Times scaled 920 

=\TimesBI scaled 860
=\TimesI scaled 860

\textfont0=\rrm  
\scriptfont0=\erm 
\scriptscriptfont0=\srm

\def\Augment#1#2{%
    \toks0\expandafter{#1}\toks2{#2}%
    \edef#1{\the\toks0\the\toks2}}

 \font\twelverma=\Times  scaled 1200
 \font\tenrma=\Times  scaled 1000
 \font\ninerma=\Times scaled 920
 =\Times scaled 840
 \font\sevenrma=\Times scaled 760
 =\Times scaled 680
 \font\fiverma=\Times scaled 600

 \Augment\tenpoint{%
  \textfont0=\tenrma  \scriptfont0=\sevenrma  
  \scriptscriptfont0=\fiverma  }

 \Augment\ninepoint{%
  \textfont0=\ninerma  \scriptfont0=\sevenrma 
  \scriptscriptfont0=\fiverma}

 \Augment\twelvepoint{%
  \textfont0=\twelverma  \scriptfont0=\ninerma  
  \scriptscriptfont0=\sevenrma}

\mathsurround=1pt
\hsize=13.45truecm
\vsize=19.5truecm
\hoffset=1.25truecm
\voffset=2truecm
\advance\baselineskip by 2pt

\predefine\til{\~}
\def\~#1{\relax\ifmmode\widetilde{#1}\else\til{#1}\fi}

\redefine \ge{\geqslant}
\define \wt#1{\mathaccent"0365{#1}}
\define \wh#1{\mathaccent"0362{#1}}

\define \iss{\,\Mathaccent{\raise -.8 ex\hbox{$\widetilde{}$\kern.1em}}\rightarrow\,}

\define \prlim{{\varprojlim}\vphantom{i}\,}
\define \inlim{{\varinjlim}\vphantom{i}\,}

\define \dimm{\operatorname{\fam0 dim\,}}

\define \im{\mathop{\fam0 im}}

\define \ab{\mathop{\fam0 ab}}

\define \degg{\operatorname{\fam0 deg\,}}

\define \res{\operatorname{\fam0 res}}

\define \Gal{\mathop{\fam0 Gal}}

\define \Pic{\mathop{\fam0 Pic}}

\define \Spec{\mathop{\fam0 Spec}}

\Mas
\HideEdStuff
\rm 
 
%%%% For GT headers and footers:

\def\issn{{\nineit ISSN 1464-8997 (on line) 1464-8989 (printed)}}

\def\gtp{{\nineit Published 10 December 2000: \ \copyright\ Geometry \& 
Topology Publications}}

\def\gtv3{{\nineit Geometry \& Topology Monographs, Volume 3 (2000) --
Invitation to higher local fields}}

%%%%% For section idents:

\def\lione
{{\rms Geometry \& Topology Monographs}}

\def \litwo{{\rms Volume 3: Invitation to higher local fields
}} 

\def\tinfo #1.#2.#3-#4
{{
\noindent  {\lione} \hfill 
\par 
\vskip-1.5pt
\noindent {\litwo} \hfill
\par 
\vskip-1,5pt
\noindent {\rms Part #1, section #2, pages #3--#4} \hfill
\vskip24pt 
}}

\def\tinfos #1.#2.#3-#4
{{
\noindent  {\lione} \hfill 
\par 
\vskip-1.5pt
\noindent {\litwo} \hfill
\par 
\vskip-1.5pt
\noindent {\rms Pages #3--#4} \hfill
\vskip24pt 
}}

\def\tinfoi #1
{{
\noindent  {\lione} \hfill 
\par 
\vskip-1.5pt
\noindent {\litwo} \hfill
\par 
\vskip-1.5pt
\noindent {\rms Pages iii--xi: Introduction and contents} \hfill
\vskip26pt 
}}

%%%% Set headers and footers %%%%

  \def\titlepagehead{\hfil}

  \newif\iftitlepage\titlepagefalse
  \newif\ifblankpage\blankpagefalse
  \def\makeheadline{
     \ifblankpage{}\else%
     \iftitlepage
\vbox{\line{\vbox to 8.5pt{}
\ninerm
\copy\HLinebox \hfill
\hglue5mm\ninebf\folio 
\titlepagehead}}%
      \else
\vbox{\ifodd\pageno\rightheadline\else\leftheadline\fi}%
      \fi\vskip 12pt\fi}%
     \def\rightheadline{\line{\vbox to 8.5pt{}%
      \ninerm
\copy\TitleBox \hfill
\hglue5mm\ninebf\folio}}%
     \def\leftheadline{\line{\vbox to 8.5pt{}%
        \unskip\ninerm\unskip\ninebf\folio\hglue5mm
      %*%
 \hfill \copy\AuthorBox
%\hfill
}}

 \footline={\ifblankpage{}\else
\iftitlepage\ninepoint\sam\hfill%} 
\line{\vbox to 8.5pt{}%\ninerm
\copy\TFLinebox
\hfill
\hglue5mm %\ninebf\folio
}
            \else
\ninepoint\sam\hfill%}
\line{\vbox to 8.5pt{}%\ninerm
\copy\FLinebox
\hfill 
\hglue5mm
}
\hfil\fi\global\titlepagefalse\fi}

\def\blankpage{{\blankpagetrue\noindent\hbox to 10pt{\hss}\vfill
\pagebreak}}

\tenpoint\rm %% always start here
 
  %%% all done and macros loaded!